\documentstyle[12pt]{article}
\makeatletter
\renewcommand{\@cite}[2]{[{{\bf #1}\if@tempswa,#2\fi}]}
\renewcommand{\@biblabel}[1]{[#1]\hfill}
\makeatother
\newtheorem{defin}{Definition}
\newtheorem{prop}{Proposition}
\newtheorem{nt}{Remark}
\newtheorem{th}{Theorem}
\newtheorem{lemma}{Lemma}
\newtheorem{cons}{Corollary}
\newcommand{\om}{\tilde{\Omega}}
\newfont{\sdbl}{msbm9}
\newfont{\dbl}{msbm10 at 12pt}
\newcommand{\eqdef}{\stackrel{\rm def}{=}}
\newcommand{\proof}{{\bf Proof.\ }}

\newcommand{\oo}{{\cal O}}
\newcommand{\ee}{{\cal E}}
\newcommand{\res}{\mathop {\rm res}}
\newcommand{\dm}{\mathop {\rm dim}}
\newcommand{\lm}{\mathop {\rm lim}}
\newcommand{\tr}{\mathop {\rm Tr}}
\newcommand{\nm}{\mathop {\rm Nm}}
\newcommand{\fdc}{\mathop { f_*^{x,C}}}
\newcommand{\fdf}{\mathop {f_*^{x,F}}}
\newcommand{\fcc}[2]{{f_*({#1},{#2})_{x,C}}}
\newcommand{\fcf}[2]{{f_*({#1},{#2})_{x,F}}}
\newcommand{\fcfn}{f_*(\;,\;)_{x,F}}
\newcommand{\fccn}{f_*(\;,\;)_{x,C}}
\newcommand{\sda}{{\mbox{\sdbl A}}}
\newcommand{\da}{{\mbox{\dbl A}}}
\newcommand{\dz}{{\mbox{\dbl Z}}}
\newcommand{\vsgm}{\varsigma}
\newcommand{\zt}{\zeta}
\newcommand{\vph}{\phi}
\newcommand{\ph}{\varphi}
\newcommand{\ps}{\psi}
\newcommand{\veps}{\varepsilon}
\newcommand{\pat}{\frac{\partial}{\partial t}}
\newcommand{\Soo}{{\rm Spec} \, \hat{\oo}_{x,X}}
\newcommand{\Ker}{{\rm Ker}\:}

\newcommand{\ka}{K_2'(\da_X)}
\newcommand{\kad}{K'_2(\da_1)}
\newcommand{\kadd}{K'_2(\da_2)}
\newcommand{\kado}{K'_2(\da_X(0))}
\newcommand{\kao}{K'_2(\oo_1)}
\newcommand{\kaoo}{K'_2(\oo_2)}

\begin{document}
\author{D. V. Osipov}
\title{Adelic constructions of direct images for differentials and
symbols\footnote{That is a modified
english version of the article appeared in "Matematicheskiy
Sbornik" 5(188) (1997).}}
\date{}
\maketitle

Let $f$  be a projective morphism from a smooth algebraic
surface $X$ to a  smooth algebraic  curve $S$ over a perfect field $k$.
Using the adelic language
we give some relative  constructions of residues  and symbols
and apply them to the Gysin morphism for differentials and
algebraic cycles.

The first section of this article
is devoted to various constructions of relative
residue maps from differentials of 2-dimensional
local fields to differentials of 1-dimensional local fields.
(See definitions~\ref{ado1}, \ref{ado2}.)
For these maps we prove some reciprocity laws.
(See propositions~\ref{pr6} and~\ref{pr7}.)
Then using the adelic resolutions of sheafs $ \Omega^2_X  $  and
$\Omega^1_S $  we apply these residue maps
to the  construction  of Gysin maps from
$H^n (X,\Omega^2_X)$ to $H^{n-1} (S,\Omega^1_S)$  ($n=1,2$).
(See propositions~\ref{prf}, \ref{prop2}.)

In the second section of this article
we assume $char k =0$
and construct relative maps $\fccn$
from $K_2$-groups of 2-dimensional local
fields $K_{x,C}$  on the surface $X$,
associated with  pairs:
an irredicuble curve $C \subset X$  and a point $x \in C$,
to  multiplicative groups of complete local fields $K_s$
of points $s$  on the curve $S$. (See theorems~\ref{t1} and~\ref{t2}.)
These  relative symbol maps are directly  connected
 to the other local maps on the surface $X$
and the  curve $S$:
such as the 2-dimensional tame symbol, 2-dimensional
residue map and so forth.
Also, we prove some relative reciprocity laws. (See corollaries from
theorem~\ref{t1}.)

Let us remark,
if $C$ is not a fibre of morphism $f$,
then the required symbol map is the usual tame symbol.
If $C$ is in the fibre of morphism $f$,
then the required symbol was originally
introduced by K.~Kato in~\cite{K}.
We give another proofs of all theorems,
we need on this symbol.
Moreover,
in theorem~\ref{t2}
we give an explicit formula for this symbol
when $char k =0$.

The third section of this article is similar to the end
of the first section.
In this part we apply the constructed symbol
maps for a construction of direct image maps
from $H^n(X, {\cal K}_2 (X))$ to $H^{n-1}(S, {\cal K}_1 (S))$ ($n=1,2$).
(See proposition~\ref{qqu}.) (Here ${\cal K}_2(X)$
(correspondingly ${\cal K}_1(S)$) is the sheaf on the surface $X$
(corr. on the curve $S$) associated to the presheaf $\{U \mapsto K_2(U) \}$
(corr. $\{U \mapsto K_1(U) \}$).)
If $n=2$, then this map is the Gysin map from   $C\!H^2(X)$
to $C\!H^1(S)$. (See proposition~\ref{ep}.)
For this goal we construct  a $K_2$-adelic
resolution of the sheaf ${\cal K}_2(X)$. (See theorem~\ref{th3}.)

Note also that all constructions in this paper
are presented by means of explicit expressions
and for almost all statements we give variants of their
proofs which don't use the higher Quillen $K$-theory.

 I would like to express my deep gratitude to my scientific
advisor, A.~N.~Parshin, for individed attention to this work.
I wish also to thank A.~B.~Zheglov for his valuable comments.

\section{Differentials and their direct images.}
\subsection{Continuous differentials.}  \label{ss1.1}
\label{cont}

 Let $K$ be a field of discrete valuation $\nu_K$.
Then by ${\cal O}_K$ denote its discrete valuation ring and by
$m_K$ denote the maximal ideal of this ring.
Let $\bar K = {\cal O}_K / m_K$ be its residue field and
$\pi : {\cal O}_K \to \bar K$ be the canonical map.

 Now if $K=k((t))$ is a local field of dimension 1,
then ${\cal O}_K = k[[t]]$, $m_K=t \cdot k[[t]] $ and $\bar K =k$.
Let $\Omega^1_{K/k}$ be the module of 1-differentials of $K$.
Put $Q=\mathop{\bigcap}\limits_{n\ge0}m_K^nd({\cal{O}}_K)$.

\begin{defin}[Continuous differentials]
$\tilde{\Omega}^1_{K/k}
\stackrel{\rm def}{=} \Omega^1_{K/k}/Q$, where $K=k((t))$.
\end{defin}

\begin{prop}
\label{p1}
Let $K=k((t))$,
then $\dm_K \tilde{\Omega}^1_{K/k} =1$,
$dt$ is a basis of $\tilde{\Omega}^1_{K/k}$ over $K$,
and for any $f \in K$ $df=\frac{\partial f}{\!\!\partial t}\,dt$.
\end{prop}
\proof
See~\cite[ ch.2.11]{S}.

 In the sequel assume that $K$ is a 2-dimensional local field,
moreover let $K=k((u))((t))$.
Then ${\cal O}_K = k((u))[[t]]$ ,
$m_K = t \cdot {\cal O}_K$,
$\bar K = k((u))$, ${\cal O}_{\bar K}=k[[u]] $.
Put $F=k((t))$, $I=k[[u,t]]$.
Then let $\varphi : I \rightarrow {\cal O}_{\bar K}$ be the quotient  map
on the ideal ($t \cdot I$).
Let  $\Omega^1_{K/k}$ be the module of 1-differentials of $K$.
Let $P_i$ be a subgroup of  $\Omega^1_{K/k}$ generated by elements
$(\varphi^{-1}(m_{\bar K}))^i d I $,
$T_j$ be a subgroup of $\Omega^1_{K/k}$  generated by elements
$       m_K^j d({\cal{O}}_K)$.
Put   $Q_I=K \cdot \mathop{\bigcap}\limits_{i,j \ge0}(P_i + T_j) $.

\begin{defin}[Continuous differentials]  \label{hh}
${_I\tilde{\Omega}^1_{K/k}}
\stackrel{\rm def}{=} \Omega^1_{K/k}/Q_I$,
where $K=k((u))((t))$
\end{defin}

The proof of the following proposition is similar to the proof
of proposition~\ref{p1}.

\begin{prop}
Let $K=k((u))((t))$,
then ${_I\tilde{\Omega}^1_{K/k}}$ is a two dimensional vector
space over the field $K$.
$du$ and $dt$ are a basis of ${_I\tilde{\Omega}^1_{K/k}}$
over $K$, and for any $f \in K$
$df=\frac{\partial f}{\!\!\partial u}\, du +
\frac{\partial f}{\!\!\partial t}\, dt$.
\end{prop}

\begin{defin} \label{hhh}
${_I\tilde{\Omega}^2_{K/k}}
\stackrel{\rm def}{=}{_I\tilde{\Omega}^1_{K/k}}\mathbin{\bigwedge_K}
{_I\tilde{\Omega}^1_{K/k}}$.
\end{defin}

\begin{nt} {\em
\label{n1}
If we choose the other representation of $K=k((u'))((t'))$,
then the subring $I$ can be changed in $K$.
Therefore our definition of ${_I\tilde{\Omega}^*_{K/k}}$
really depends on the choice of an embedding $I \hookrightarrow K$.
Moreover, if $\phi_I$ is the natural homomorphism from
${\Omega}^*_{K/k}$ to ${_I\tilde{\Omega}^*_{K/k}}$,
then the composition of $\phi_I$ with the residue map
(see defenitions 4 and 5) will depend on the choice of an embedding
$I \hookrightarrow K$.

But the indeterminancy of such kind can be removed when the two-dimensional
local field $K$ appears from geometrical datas. In this case
there is a canonical embedding of the ring $I$ to the field $K$.
(See section~\ref{ss1.2}).
} \end{nt}

\begin{nt} {\em
There is an other approach to the definition of continuous
differentials: by introducing  some topology on the module ${\Omega}^*_{K/k}$
A. Yekutieli  considered in~\cite{Y} the module
of the "right" differentials  ${\Omega}^{*,sep}_{K/k}$
for the multidimensional local fields.
One can prove, that our ${_I\tilde{\Omega}^*_{K/k}}$
coinside with  ${\Omega}^{*,sep}_{K/k}$.
(The topology of the field $K$ is the following: the set of
subgroups
 $J_{i,j}=
(\varphi^{-1}(m_{\bar K}))^i I + t^j {\cal O}_K$ is
a basis of neighbourhourds of $0$ in the additive topology of the ring
${\cal O}_K$,
and $K= \lm\limits_{n \to} \; t^{-n} \cdot {\cal O}_K $
is equipped now by the topology
of inductive limit. See also~\cite{PF}.)
} \end{nt}

\begin{defin}[A residue map]  \label{ado1}
Let $\omega \in
{_I\tilde{\Omega}^2_{K/k}}$, $\omega=\sum\limits_i \omega_i (u)\wedge
t^i dt$, then $\res_{K/\bar K} \omega \eqdef \mathop{\rm res}_t \omega
\stackrel{\rm def}{=} \omega_{-1}(u) \in \tilde{\Omega}^1_{\bar K}$.
\end{defin}

The map $\res_{K/{\bar K}}$ is well defined due to the following
proposition.

\begin{prop}
If we fix an embedding $\bar K \hookrightarrow K$,
then $\mathop{\rm res}_t $ does not depend on the choice
of the local parameter $t \in I$.
\end{prop}
\proof
Actually this fact is about 1-dimensional local fields and follows
from~\cite[ ch.~2.]{S}).

\begin{defin}[A residue map] \label{ado2}
Let $\omega \in {_I\tilde{\Omega}^2_{K/k}}$,
 $\omega=\sum\limits_i u^i du\wedge
\omega_i(t)$,
then $\res_{K/F} \eqdef \mathop{\rm res}_u \omega
\stackrel{\rm def}{=} \omega_{-1}(t) \in \tilde{\Omega}^1_F$.
\end{defin}

The map $\res_{K/F}$ is well defined due to the following proposition:

\begin{prop}
If we fix an embedding $F \hookrightarrow K$,
then $\mathop{\rm res}_u $ does not depend on the choice of the local
parameter $u \in I$.
\end{prop}
\proof
Let $u'$ be an other local parameter, i.~e. $u'\in I$,
$\nu_{\bar K}(\pi(u')) = 1$.
Let $$\mathop{\rm res}\nolimits_u \omega=\sum_i a_i t^i dt
\qquad \quad
\mbox{and}
\qquad \quad
\mathop{\rm res}\nolimits_{u'}\omega=\sum_i a_i' t^i dt {\mbox .}$$
Then we have
$$ a_i=\mathop{\rm res}\nolimits_{F/k}(t^{-1-i}
\mathop{\rm res}\nolimits_u(\omega))=
\mathop{\rm res}\nolimits_{F/k} \cdot \mathop{\rm res}\nolimits_u
(t^{-1-i}\omega)=
\mathop{\rm res}\nolimits_{t,u}(t^{-1-i}\omega)$$
$$ a_i'=\mathop{\rm res}\nolimits_{F/k}(t^{-1-i}
\mathop{\rm res}\nolimits_{u'}(\omega))=
\mathop{\rm res}\nolimits_{F/k} \cdot \mathop{\rm res}\nolimits_{u'} (t^{-1-i}\omega)=
\mathop{\rm res}\nolimits_{t,u'}(t^{-1-i}\omega) {\mbox ,}$$
where $\mathop{\rm res}_{F/k} : \tilde{\Omega}^1_{F/k} \to k$
is an usual residue map for the 1-dimensional local field described
in~\cite[ ch.~2]{S}), $\mathop{\rm res}_{t,u}$ and
$\mathop{\rm res}_{t,u'}$ $:\tilde{\Omega}^2_{K/k} \to k$
are residue maps of 2-dimensional local fields:
$\mathop{\rm res}_{t,u}(\sum a_{i,j}u^it^j du
\wedge dt)= a_{-1,-1}$. But from~\cite[ ch. 2.5.3]{PF}
(or \cite{P2}) $\res_{K/k} \eqdef \res_{t,u}$ does not
depend on the choice of local parameters, i.~e.
$\res_{t,u}=\mathop{\rm res}_{t,u'} $.
Therefore $a_i=a_i'$ and $\res_u = \res_{u'}$.
The proof is finished.

\subsection{Direct images.} \label{ss1.2}
Let $k$ be a perfect field,
$X$ be a smooth algebraic surface over $k$,
$S$ be a smooth algebraic curve over $k$,
$f:X \to S$ be a projective morphism over $k$.
In the sequel assume $f:X \to S$ is a smooth morphism with
connected fibres, but remark that the all further constructions
are transfered to the nonsmooth case without difficulties.

Now let $\hat {\cal O}_s$  be a complete local discrete valuation
ring at point  $s \in S$ ,
let  $K_s$  be its fraction field.
For the sake of the simplicity of further notations suppose
$k(s)=k$.
Let  $\tau$  be a local parameter at point $s$,
then $\hat {\cal O}_s=k[[\tau]]$ and $K_s=k((\tau))$.
Now for any closed point $x \in F$,
where $F$ is the fibre of $f$: $F=f^{-1}(s)$,
and any irreducible curve $C \in X$  such that $x \in C$ and
$x$ is a nonsingular point of $C$ (or $C$ has only one
analytic branch at $x$)  one can construct a canonical
2-dimensional local field $K_{x,C}$. (See~\cite{PF} or
\cite{P2}, \cite{L}). Let us consider two cases.

\begin{enumerate}
\item      If $C\ne F$, then  $K_{x,C} \simeq k(C)_x((t_C))$,
where $t_C=0$ is a local equation of the curve $C$ at the point $x$,
$k(C)_x$  is the completion of the function field of
the curve $C$ at the point $x$.
Besides, $k(C)_x$ is a finite extension of  $K_s$  under the map $f^*$.

\item  If       $C=F$, then  $K_{x,F} \simeq k'((u))((t)) $,
where $k'=k(x)$,
$t=f^*(\tau)$ and $u$ is from $k(X)$ such that $u$ and $t$
are local parameters at the point $x$.
Besides, $k'((t))$ is a finite separable extension of $K_s$
under the map $f^*$.
\end{enumerate}

Now let ${\cal O}_{x,X}$ be a local ring of the point $x$
on surface $X$, let
$\hat {\cal O}_{x,X}$
be a completion of the  ring  ${\cal O}_{x,X}$
at maximal ideal. Then for any irredicuble curve  $C \ni x$  there exists
a canonical embedding of the ring $\hat {\cal O}_{x,X}$
to the field $K_{x,C}$.
But the ring $\hat {\cal O}_{x,X}=k'[[u,t]]$,
and it is exactly the ring $I$ from section~\ref{ss1.1}.
Therefore in the sequel we shall suppose $I=\hat {\cal O}_{x,X}$
and we shall write $\tilde{\Omega}^*$ instead of ${_I\tilde{\Omega}^*}$.

\begin{defin}[Direct image maps $\fdc : \om^2_{K_{x,C}/k'} \to \om^1_{K_s/k}$] \label{d6}

\begin{enumerate}

\item If $C \ne F$, then for any $\omega \in \om^2_{K_{x,C}/k'}$
$$ \fdc(\omega)\eqdef \tr\nolimits_{k(C)_x/K_s}
\cdot \mathop{\rm res}\nolimits_{K_{x,C}/k(C)_x}(\omega)
\in \om^1_{K_s/k} \mbox{.}$$
\item If $C=F$, then for any $\omega \in \om^2_{K_{x,F}/k'}$
$$ \fdf(\omega)\eqdef \tr\nolimits_{k'((t))/K_s}
\cdot \mathop{\rm res}\nolimits_{K_{x,F}/k'((t))}(\omega)
\in \om^1_{K_s/k} \mbox{.}$$
\end{enumerate}
\end{defin}
\begin{nt} {\em
Here $\tr$ is a trace map of differential forms.
See its definition in~\cite[ ch.~2]{S}, also in~\cite{L} and~\cite{PF}
for  inseparable case.
} \end{nt}

If $x$ is a singular point of the irredicuble curve $C$,
then  $$  K_{x,C} \eqdef        \bigoplus_i K_{x,C_i}  $$
			$$  \om^2_{K_{x,C}} \eqdef \bigoplus_i  \om^2_{K_{x,C_i}}  $$
 and we put                     $$
      f_*^{x,C} \eqdef \sum_i f_*^{x,C_i}  \mbox{,}$$
where $C_i$ denote the "analytic" branches of $C$ at $x$.

From definition~\ref{d6} we obtain easily the following proposition:
\begin{prop} \label{pr5}
Let $\omega \in \om^2_{\oo_{K_{x,C}}/k'} \subset \om^2_{K_{x,C}/k'}$,
i.~e.  $\omega = g \: dt_1 \wedge dt_2$
$\nu_{K_{x,C}}(g) \ge 0$  for $K_{x,C}=k'((t_1))((t_2))$
and $t_1, t_2 \in \hat {\cal O}_{x,X}$; then

\begin{enumerate}
\item If $C \ne F$, then $\fdc(\omega)=0$.
\item If $C = F$, then $\fdf(\omega) \in \om^1_{\oo_{K_s}/ k}$, i.~e.
$\fdf(\omega)=\varepsilon $, where $\varepsilon=h d\tau$ and
$\nu_{K_s}(h) \ge 0$.
\end{enumerate}
\end{prop}

For any irredicuble curve $C \in X$ let $K_C$  be the completion
of the field $k(X)$ with respect to the discrete valuation defined
by $C$, i.~e. $K_C \eqdef k(C)((t_C)) \subset K_{x,C}$.
By $\oo_C$ denote the ring $\oo_{K_C}$, i.~e. $\oo_C \eqdef k(C)[[t_C]]$.
Note also  that for any point $x \in C$ we have a natural map
from $\Omega^*_{K_C/k}$ to $\om^*_{K_{x,C}/k'}$.

\begin{prop}[The reciprocity law along a fibre.]   \label{pr6}
Let $\omega \in \Omega^2_{K_F/k}$, then the following series converges
in the topology defined by the discrete valuation of the point $s$
and
\begin{equation}
\label{sum}
\sum_{x \in F} \fdf (\omega)=0 {\mbox .}
\end{equation}
\end{prop}
\begin{nt} {\em \label{n4}
It can be that the sum~(\ref{sum}) is really infinite.
That is an example (see also~\cite{K}).
\begin{tabbing}
Let \qquad\=$X={\bf P}^1 \times {\bf A}^1$ and $(u,t)$ are their
coordinates,\\
            \>$S={\bf A}^1$ with the coordinate $t$,\\
            \>$F={\bf P}^1$ and the point $s=(0,0)$.
\end{tabbing}
Then for $\omega=(u^{-1}+\frac1{u-1}t+\frac1{u-2}t^2+
\frac1{u-3}t^3+\ldots)du\wedge dt \in \om^2_{K_F/k}$,
\begin{tabbing}
\= at point $x=(l,0)$ \quad  \= we have \qquad \= $ \fdf(\omega)=t^ldt $,
where $l=0,1,2,\ldots $ \\
\> at point $x=(\infty,0)$ \> we have  \> $ \fdf(\omega)=(-1-t-t^2-t^3-\ldots)dt$\\
\> at other points $x$ \> we have     \> $\fdf(\omega)=0  \mbox{.}$
\end{tabbing}
} \end{nt}
\proof[of proposition \ref{pr6}]

Let $n=\nu_{K_F}(\omega)$, i.~e. if $\omega=g\: du\wedge dt$  in
some  $\om^2_{K_{x,F}/k'}$, then $n \eqdef \nu_{K_{x,F}}(g)$
and this number does not depend on the choice of the point $x \in F$.

\begin{enumerate}
\item \label{b1}
It is clear, that for any point $x \in F$ if $\fdf(\omega)=hdt$, then
$\nu_{K_s}(h) \ge n$.
\item  \label{b2}
For any $m \ge n$ there exists only a finite number of
points  ${x \in F}$ such that if
$\fdf(\omega)=hdt$ $h=a_nt^n+a_{n+1}t^{n+1}
+\ldots+a_mt^m+\ldots$, then   $a_m \ne 0$. And the sum over all
such $a_m$ is equal to $0$ for each fixed $m$.

Indeed, if we fix the point $x \in F$ and $m \ge n$, then
$$
\begin{array}{rcl}
a_m &=&\mathop{\rm res}_{k((t))/k}(t^{-1-m} \cdot
\fdf(\omega))=\\
&=&\mathop{\rm res}_{k((t))/k}(t^{-1-m} \cdot
\tr_{k'((t))/k((t))}\mathop{\rm res}_{K_{x,F}/k'((t))}(\omega))=\\
&=&\mathop{\rm res}_{k((t))/k}(\tr_{k'((t))/k((t))}\mathop{\rm res}_{K_{x,F}/k'((t))}(t^{-1-m}
\cdot \omega))=\\
&=&\mbox{(see remark~\ref{n5})}=\\
&=& \tr_{k'/k}\mathop{\rm res}_{k'((t))/k'}
\mathop{\rm res}_{K_{x,F}/k'((t))}(t^{-1-m} \cdot \omega)=\\
&=&\mathop{\rm res}_{x,F}(t^{-1-m}
\cdot \omega) \mbox{.}
\end{array}
$$
Note that $(t^{-1-m} \omega) \in \Omega_{K_F/k}$, therefore the
statement~\ref{b2}
follows from the reciprocity law  along a projective curve for
the 2-dimensional residue map:
for any $\eta \in \Omega^2_{K_F/k}$  $\sum\limits_{x \in F} \res_{K_{x,F}/k}=0$
and the number of terms $\ne 0$ in this sum is finite.
(See~\cite[ ch. 4.1]{PF}, \cite{P2}, \cite{L}).
\end{enumerate}

Now proposition~\ref{pr6} follows from statements~\ref{b1} and \ref{b2}.

\begin{nt} {\em \label{n5}
In this proof we used the commutativity of the operations
 $\mathop{\rm res}$ and $\tr$.
(See~\cite[ ch.~2]{S}).
} \end{nt}

For any point $x \in X$ let $K_x \eqdef \hat \oo_{x,X}
\cdot
k(X)$ be the subring in fraction field of $\hat \oo_{x,X}$.
For any irredicuble curve $C \ni x$ there exist canonical maps
from $K_x$ to $K_{x,C}$ and from $\Omega^*_{K_x/k}$
to $\om^*_{K_{x,C}/k}$.

\begin{prop}[The reciprocity law around a point]  \label{pr7}
Fix a point ${x \in F}$, $\omega \in \Omega^2_{K_x/k}$.
Then in the following infinite series the number of
terms $\ne 0 $ is finite and

\begin{equation} \label{sum7}
\sum_{C\ni x} \fdc(\omega)=0 \mbox{.}
\end{equation}

(The last sum is over all irredicuble curves $C \subset X$
such that $x \in C$.)

\end{prop}
\proof

It is clear, that if $C \not \subset Supp(\omega)$, then $\fdc(\omega)=0$.
Hence in the sum~(\ref{sum7}) the number of terms $\ne 0$
is finite. (Compare with proposition~\ref{pr6}).

Let $\sum\limits_{C \ni x} \fdc(\omega)=(\sum \limits_m a_mt^m)dt$, then
$$
\begin{array}{rcl}
a_m &=& \mathop{\rm res}_{k((t))/k}(t^{-1-m} \cdot \sum \limits_{C \ni x}
\fdc(\omega))=\\
&=&\mathop{\rm res}_{k((t))/k}(t^{-1-m} \cdot (\sum\limits_{C \ni x}
\tr_{k(C)_x/k((t))}  \mathop{\rm res}_{K_{x,C}/k(C)_x}(\omega)))=\\
&=&\mathop{\rm res}_{k((t))/k}(\sum\limits_{C \ni x}  \tr_{k(C)_x/k((t))}
t^{-1-m} \cdot \mathop{\rm res}_{K_{x,C}/k(C)_x}(\omega))=\\
&=& res_{k((t))/k}(\sum\limits_{C \ni x}   \tr_{k(C)_x/k((t))}
 t^{-1-m} \cdot \mathop{\rm res}_{K_{x,C}/k(C)_x}(\omega))=\\
 &=& \sum\limits_{C \ni x} \tr_{k'/k} \mathop{\rm res}_{k(C)_x/k'}
 \mathop{\rm res}_{K_{x,C}/k(C)_x}(t^{-1-m} \cdot \omega)=\\
 &=&\sum\limits_{C \ni x} \mathop{\rm res}_{x,C}(t^{-1-m} \cdot \omega) =\\
 &=& 0
 \\
\end{array}
$$

The last follows from the reciprocity law around the point for
2-dimensional residue map:
for any $\eta \in \Omega^2_{K_x/k}$  $\sum\limits_{C \ni x} \res_{K_{x,C}/k}=0$.
(See~\cite[ ch. 4.1]{PF}, \cite{P2}).

\subsection{Adelic differentials and the Gysin morphism.}

Retain all notations of the previous section.

\begin{defin}
Let $S$ be a smooth curve, then
  $$
   \Omega^1_{{\mbox{\sdbl A}}_S} \eqdef
   \left \{ (\ldots , f_s d\tau_s, \ldots)
   \in \prod_{s \in S} \om^1_{K_s/k}\;\mbox{, where \quad}
   \nu_{K_s}(f_s) \ge 0 \mbox{ almost everywhere} \right \} \mbox{.}
  $$
  Here $\tau_s$ is a local parameter
  at the point $s$ and $\nu_s$ is the corresponding
  discrete valuation.

  For any divisor $D$ of $S$ suppose
  $$
  \Omega^1_{\sda_S}(D) \eqdef \{f \in \Omega^1_{\sda_S} \;
  \mbox{, where } \nu_s(f_s) \ge - \nu_s(D) \} \mbox{.}
  $$
\end{defin}
\begin{defin}     \label{d8}
Let $X$ be a smooth surface and a pair $x \in C$ runs over all
irredicuble curves $C \subset X$ and all points $x \in C$, then
$$
\Omega^2_{\sda_X} \eqdef
\left \{ (\ldots,\omega_{x,C},\ldots) \in \prod_{x \in C} \om^2_{K_{x,C}/k}
\; \mbox{ under the following conditions} \right\}
$$
If $ \omega_{x,C}=\sum\limits_{i \ge \nu(\omega_{x,C})}
\omega_i^{x,C}(u_{x,C}) \wedge t^i_C dt_C$,
where $t_C=0$ is local equation of the curve $C$ and for all points $x$
from some open smooth  set $U$ of the curve $C$
$$K_{x,c}=k(x)((u_{x,C}))((t_C)) \qquad , \qquad
 \omega_i^{x,C}(u_{x,C}) \in \om^1_{k(x)((u_{x,C}))} \mbox{.}$$
Then
\begin{enumerate}
\item there exists a divisor $D=\sum \nu_C(D) \cdot C$ on $X$
such that $\nu(\omega_{x,C}) \ge \nu_C(D)$.
\item For  fixed curve $C \subset X$ and  $i$
the collection $\omega_i^{x,C}(u_{x,C})$ belongs to $\Omega^1_{\sda_U}$
when $x$ runs $U$.
\end{enumerate}
\end{defin}

\begin{defin}
Let $D$ be a divisor on $X$, then
$$
\Omega^2_{\sda_X}(D) \eqdef \{ \omega \in \Omega^2_{\sda_X}\;
\mbox{, where} \quad
\nu(\omega_{x,C}) \ge -\nu_C(D)\} \mbox{.}
$$
\end{defin}

Now {\em define} the map      \label{adf}
$f_* \eqdef \mathop{\sum\limits_{f(x)=s}}\limits_{C \ni x}\fdc $
from  $\Omega^2_{\sda_X} \subset \prod\limits_{x \in C} \om^2_{K_{x,C}/k'}$
to $\Omega^1_{\sda_S} \subset \prod\limits_{s \in S} \om^1_{K_s/k}$.

\begin{prop}
The map  $f_* \eqdef \mathop{\sum\limits_{f(x)=s}}\limits_{C \ni x}\fdc $
from $\Omega^2_{\sda_X}$ to $\Omega^1_{\sda_S}$ is
well defined, i.~e. this infinite series converges
at every point $s \in S$.
\end{prop}
\proof

Over each point $s \in S$ we have
\begin{equation}     \label{bsum}
 \mathop{\sum\limits_{f(x)=s}}\limits_{C \ni x}\fdc =
   \mathop{\sum\limits_{f(x)=s}}\fdf +
  \mathop{ \mathop{\sum\limits_{f(x)=s}}\limits_{ C \ni x}}\limits_{C \ne F}
  \fdc  \quad \mbox{.}
 \end{equation}

Now due to adelic conditions on elements of $\Omega^2_{\sda_X}$
the first sum from the right part of expression~(\ref{bsum})
converges. That can be proved by the same method as proposition~\ref{pr6}.

The second sum contains only a finite number terms $\ne 0$.
It follows easily from the definitions of $\fdc$ (when $C \ne F$)
and of $\Omega^2_{\sda_X}$.

\begin{nt} {\em
Notice that the expression
$\mathop{\sum\limits_{f(x)=s}}\limits_{C \ni x}\fdc $
applied to the whole $\prod\limits_{x \in C} \om^2_{K_{x,C}/k'}$
make not sense, since this series will not converge.}
\end{nt}

For the curve $S$ consider the following complex $\Omega^1({\cal A}_S)$:

$$
\begin{array}{ccc}
\Omega^1_{k(S)/k} \oplus \Omega^1_{\sda_S}(0) & \longrightarrow &
\Omega^1_{\sda_S}  \\[2pt]
(f_0,f_1) & \longmapsto & f_0 + f_1 \quad \mbox{.} \\[3pt]
\end{array}
$$

Then from~\cite{S} we have that
\begin{equation}  \label{com1}
H^*(\Omega^1({\cal A}_S)) \simeq
H^*(S,\Omega^1_S)    \mbox{,}
\end{equation}
where $\Omega^1_S$ is the sheaf of regular 1-differentials on
the curve $S$.

Using the diagonal map of $\Omega^2_{K_C/k}$ to $\prod\limits_{x \in C}
\om^2_{K_{x,C}/k}$ and of $\Omega^2_{K_x/k}$ to $\prod\limits_{C \ni x}
\om^2_{K_{x,C}/k}$ we put
$$
\Omega^2_{\sda_1} \eqdef (\prod\limits_{C \subset X} \Omega^2_{K_C/k})
 \cap \Omega^2_{\sda_X}
\quad \mbox{,} \quad
\Omega^2_{\sda_2} \eqdef (\prod\limits_{x \in X} \Omega^2_{K_x/k})
\cap \Omega^2_{\sda_X}
         $$
$$
\Omega^2_{{\cal O}_1} \eqdef (\prod_{C \subset X} \Omega^2_{{\cal O}_C}/k)
\cap \Omega^2_{\sda_X}  \; \subset \Omega^2_{\sda_1}
\quad \mbox{,} \quad
\Omega^2_{{\cal O}_2} \eqdef  (\prod_{x \in X} \Omega^2_{\hat{\cal O}_{x,X}/k})
\cap \Omega^2_{\sda_X} \; \subset  \Omega^2_{\sda_2}
$$

For the surface $X$ consider the following complex $\Omega^2({\cal A}_X)$:
{
\arraycolsep=2pt
$$
\begin{array}{@{}ccccc@{}}
\Omega^2_{k(X)/k} \oplus \Omega^2_{{\cal O}_1} \oplus
\Omega^2_{{\cal O}_2}
 & \longrightarrow  &
\Omega^2_{\sda_2} \oplus \Omega^2_{\sda_1}
 \oplus \Omega^2_{\sda_X}(0)  & \longrightarrow & \Omega^2_{\sda_X} \\[7pt]
(f_0,f_1,f_2) & \mapsto  &
(f_2 - f_0, f_0 + f_1, - f_1 - f_2) \\
&& (g_1,g_2,g_3) & \mapsto & g_1 + g_2 + g_3
\end{array}
$$
}

From~\cite[ ch. 4.2]{PF} (or \cite{P2}) we have
\begin{equation} \label{com2}
H^*(\Omega^2({\cal A}_X)) \simeq
H^*(X,\Omega^2_X)    \mbox{,}
\end{equation}
where $\Omega^2_X$ is the sheaf of regular 2-differentials on
the surface $X$.

Now extend the map $f_*$ to the complex $\Omega^2({\cal A}_X)$.
We have the following proposition.
\begin{prop}  \label{prf}
$f_*$ maps the complex $\Omega^2({\cal A}_X)$
to the complex
$$
0 \longrightarrow \Omega^1_{k(S)/k} \oplus \Omega^1_{\sda_S}(0)
\longrightarrow \Omega^1_{\sda_S}
$$
and this map is a morphism of complexes.
\end{prop}
{\bf Corollary} \\
{\em
$f_*$ gives us the maps from $H^1(X, \Omega^2_X)$
to $H^0(S, \Omega^1_S)$
and from $H^2(X, \Omega^2_X)$
to $H^1(S, \Omega^1_S)$ }\\[5pt]
{\bf Proof (of corollary).}
It follows from~\ref{com1} and~\ref{com2}.\\[5pt]
{\bf Proof (of proposition~\ref{prf}).}
It is enough to prove the following statements:
\begin{enumerate}
\item \label{st1}
$f_*(\Omega^2_{\sda_2})=0$.
Indeed, by proposition~\ref{pr7}
$\sum\limits_{C \ni x}\fdc=0$ for each point $x \in X$,
therefore
$\mathop{\sum\limits_{f(x)=s}}\limits_{C \ni x}\fdc=
\sum\limits_{f(x)=s}(\;\sum\limits_{C \ni x}\fdc\;) \:$.
\item  \label{st2}
$f_*(\Omega^2_{\sda_1}) \subset \Omega^1_{k(S)/k}$.
Indeed, after an application of the first sum from
the right part of expression~(\ref{bsum}) to $\Omega^2_{\sda1}$
we obtain $0$ by proposition~\ref{pr6}.
The application of the second sum from the right part
of expression~(\ref{bsum}) to $\Omega^2_{\sda1}$ yields us
elements of $\Omega^1_{k(S)/k}$.
\item
$f_*(\Omega^2_{\sda_X(0)}) \subset \Omega^1_{\sda_S(0)}$
by proposition~\ref{pr5}.
\item
$f_*(\Omega^2_{k(X)/k})=0$  and $f_*(\Omega^2_{{\cal O}_2}) =0$
by statement~\ref{st1}.\\
$f_*(\Omega^2_{{\cal O}_1}) =0$ by trivial reasons like statement~\ref{st2}.
\end{enumerate}

\begin{prop}                \label{prop2}
If $X$ and $S$ are projective varieties,
then the constructed map $f_*$ from $H^2(X, \Omega^2_X)$
to $H^1(S,\Omega^1_S)$ is dual to the
pull-back map $f^*$ from $H^0(S,{\cal O}_S)$ to $H^0(X,{\cal O}_X)$,
i.~e. $f_*$ is the Gysin map.
\end{prop}
{\bf Proof} is an easy consequence of the definition $f_*$,
an equality $\res_{K_s/k} f_*^{x,C}=\res_{K_{x,C}/k}$,
and of the following commutative diagram
$$
\begin{array}{ccccc}
k &
\stackrel{\sum \mathop{\rm res}_{K_{x,C}/k}}{\longleftarrow} &
\qquad
H^0(X,\oo_X) \simeq k & \times &
 \Omega^2_{\sda_X} \; / \;
(\Omega^2_{\sda_1}+\Omega^2_{\sda_2}+\Omega^2_{\sda_X}(0) )  \\
 \begin{picture}(2,46)
 \put(0,43){\line(0,-1){40}}
 \put(2,43){\line(0,-1){40}}
 \end{picture} &
 \qquad &
 \begin{picture}(0,46)
 \put(0,3){\vector(0,1){40}}
 \put(0,23){$\; {\rm id}$}
 \end{picture} &&
 \begin{picture}(0,46)
 \put(0,43){\vector(0,-1){40}}
 \put(0,23){$\; f_*$}
 \end{picture}\\
 k&
\stackrel{\sum \mathop{\rm res}_{K_s/k}}{\longleftarrow} &
\qquad
H^0(S,\oo_S) \simeq k & \times &
\Omega^1_{\sda_S} \; /\; (\Omega^1_{k(S)/k}+
\Omega^1_{\sda_S}(0))  \mbox{.}
\end{array}
$$

Here $\times$ is the cup-product;
$\sum\limits_{x \in C} \mathop{\rm res}_{K_{x,C}/k}$
gives us the Serre duality between $H^2(X,\Omega^2_X)$
and $H^0(X,{\cal O}_X)$;
$\sum\limits_s \res_{K_s/k}$ is the Serre duality
between $H^1(S, \Omega^1_S)$ and $H^0(S,{\cal O}_S)$.
(See~\cite[ ch. 4.3]{PF}.)

\section{Direct images and symbols.} \label{s2}
During this section assume $char k =0$.
All other assumptions and notations retain from section~\ref{ss1.2}.
Let us remark the following facts.

If $K=k((t))$ is 1-dimensional local field,
then
\begin{equation} \label{dec1}
K^* \simeq \{t^m\} \times {\cal U}_K \simeq \dz \times {\cal U}_K\qquad
\mbox{and} \qquad  {\cal U}_K \simeq k^* \times {\cal U}^1_K  {\mbox ,}
\end{equation}
where the group of 1-units ${\cal U}^1_K=1+m_{K}$,
i.~e. ${\cal U}^1_K
 = 1+t \cdot k[[t]]$.

If  $K=k((u))((t))$  is 2-dimensional local field,
then
\begin{equation} \label{dec2}
K^* \simeq \{t^m u^n\} \times {\cal E}_K \simeq \dz \times
\dz \times {\cal E}_K\qquad
\mbox{and} \qquad  {\cal E}_K \simeq k^* \times {\cal E}^1_K  {\mbox ,}
\end{equation}
where the group of 1-units
 ${\cal E}^1_K = 1 + \pi^{-1}(m_{\bar K})$,
i.~e.            ${\cal E}^1_K =
1+u \cdot
k[[u]] +t \cdot k((u))[[t]]$.
(Remind that $\pi : {\cal O}_K \to {\bar K}$ is the canonical map.)

If $K$ is a field, $\nu$ is the discrete valuation of $K$,
$m_K$ is the maximal ideal of the valuation ring,
$\bar K$ is the residue field of $K$,
then for any $\ph$ and $\ps$  from $K^*$
the tame symbol $(\;,\;)_K$ from $K^* \times K^*$
to $\bar{K}^*$ is
\begin{equation} \label{tc1}
  (\ph,  \ps)_K   \eqdef   ((-1)^{\nu(\ph) \nu(\ps)} \ph^{ \nu(\ps)}
  \ps^{ -\nu(\ph)}) \; \bmod m_K  \;              \quad \mbox{.}
\end{equation}

If $K=k((u))((t))$ is a 2-dimensional local field,
then for any  $\vph$, $\ph$ and $\ps$ from $K^*$
one can define an analogous symbol $(\;,\;,\;)_K$
from  $K^* \times K^* \times K^*$  to $k^*$:
\begin{equation}   \label{tc2}
(\vph,\ph,\ps)_K \eqdef (\partial_2 \: \partial_3 \:(\vph,\ph,\ps))^{-1}
\mbox{.}
\end{equation}
Here $(\vph,\ph,\ps)$ belongs to the Milnor $K$-group $K_3^M(K)$,
 $\partial_m  :  K^M_m(L) \to K^M_{m-1}( \bar L)$
is the boundary map for any local field $L$
and its residue field $\bar L$. (See~\cite{MI} or~\cite[ ch. 3.1]{PF}.)
As above, there exists an explicit formula for  $(\;,\;,\;)_K$.
(See~\cite[ ch. 3.2]{PF} or \cite[ \S 3]{P3}.)

Now we shall give the following theorem.

\begin{th} \label{t1}
Fix a point $x \in
X$, an irredicuble curve $C \ni x$ such that $x$ is
a nonsingular point of $C$  , and a point
$s \in S$ such that $f(x)=s$; then there exists a map  $\fccn$
from $K_{x,C}^*  \times  K_{x,C}^*$
to $K_s^*$ such that the following conditions hold:
\begin{enumerate}
\item for any $\ph , \ps \in K_{x,C}^* \;,\; \xi \in K_s^*$
\begin{equation} \label{d}
\nm\nolimits_{k(x)/k}(\ph , \ps , f^*(\xi ))_{K_{x,C}} =
(\fcc{\ph}{\ps}, \xi)_{K_s}
\end{equation}
\item for any $\ph , \ps \in K_{x,C}^* \; , \;\zt \in K_s$
\begin{equation} \label{dd}
\tr\nolimits_{k(x)/k}(\ph,\ps;f^*(\zt)]_{K_{x,C}} =
(\fcc{\ph}{\ps};\zt]_{K_s}
\qquad \mbox{,}
\end{equation}
where
$$
(\ph,\psi; \vsgm]_{K{x,C}}=\mathop{\rm res}\nolimits_{K_{x,C}/k}\left(\vsgm
\frac{d\ph}{\ph} \wedge \frac{d\ps}{\ps}\right) \quad
\mbox{ ,  $\vsgm \in K_{x,C}$}
$$
and
$$
(\xi;\zeta]_{K_s}=\mathop{\rm res}\nolimits_{K_s/k}\left(\zeta \frac{d\xi}{\xi}\right)
\quad
\mbox{ ,  $\xi \in K_s^* $.}
$$
\end{enumerate}
\end{th}
\begin{nt} \em If $x$ is a singular point of an irredicuble curve $C$,
then we put $\fccn = \prod\limits_i f_*(\;,\;)_{x,C_i}$,
where $C_i$ are the "analitic" branches of the curve $C$ at the point $x$.
(Compare with section~\ref{ss1.2}.)
\end{nt}

This theorem will be proved later by means of an explicit
formula for $\fccn$. (See theorem~\ref{t2}.)

\begin{cons}  \label{c1}
The map $\fccn$ is uniquely determined by the conditions
of theorem~\ref{t1}.
\end{cons}
{\bf Proof} (of corollary~1).

Let there exists two maps satisfied the conditions of theorem~\ref{t1}.
Then if we divide the first map by the second map
we obtain that there exists an element                 ${\zt \in K_s^*} \; ,
	\; {\zt \ne 1}$  such that for any element    $\xi \in K_s^*$:
\begin{equation} \label{e}
 (\zt, \xi)_{K_s}=1  \qquad   \mbox{and}  \qquad  (\zt ; \xi]_{K_s}=0
\mbox{.}
\end{equation}
But this is not right, since

\begin{enumerate}
\item If $\nu_{K_s}(\zt) \ne 0$, then for $\xi \in k^*$,
$ \xi^{\nu_{K_s}(\zt) } \ne 1$
we have $(\zt,\xi)_{K_s} \ne 1$. This contradicts expressions~(\ref{e}).
\item If $\nu_{K_s}(\zt)=0$ but $\zt \notin {\cal U}^1_{K_s}$, then
$(\zt, \tau)_{K_s} \ne 1$. (Remind $K_s \simeq  k((\tau))$).
\item
If $\zt \in {\cal U}^1_{K_s}$ and $\nu_{K_s}(\zt -1) = n$,
then $(\zt; \tau^{-n}]_{K_s} \ne 0$.
\end{enumerate}
\begin{cons} \label{cor2}
The map $\fccn$ is a symbol map, i.~e. $\fccn$
is well defined as the map:
$$
\fccn : K_2(K_{x,C}) \longrightarrow K_1(K_s) \mbox{.}
$$
\end{cons}
{\bf Proof} (of corollary~\ref{cor2}).

Note that from definitions of $(\;,\;,\;)_{K_{x,C}}$
and $(\;,\; ;\;]_{K_{x,C}}$  (or from~\cite{PF})
we obtain that the map  $(\;,\;,\;)_{K_{x,C}} :
K_{x,C}^* \times K_{x,C}^* \times K_{x,C}^* $ to $k^*$
is well defined as a map from $K_2(K_{x,C}) \times K_{x,C}^*$ to $k*$
and the map from  $(\;,\; ;\;]_{K_{x,C}} :
K_{x,C}^* \times K_{x,C}^* \times K_{x,C}$ to $k$ is
well defined as a map from $K_2(K_{x,C}) \times K_{x,C}$ to $k$.
Now the proof is done by means of the same methods as the proof of
corollary~\ref{c1}.

\begin{cons}[The reciprocity law along a fibre.]     \label{c3}
Fix a point $s \in S$,
the fibre $F = f^{-1}(s)$,
$\ph, \ps \in K_F^*$.
Then the following infinite product converges in $K_s^*$ and
$$ \prod_{x \in F} \fcf{\ph}{\ps}=1  \mbox{.}$$
 \end{cons}
{\bf Proof} (of corollary~\ref{c3}).

From~\cite[ ch. 7.1]{PF} (or~\cite{UMN2} and~\cite{P2})  we have analogous reciprocity laws
along a projective curve $F$ on the surface $X$
for $(\;,\;,\;)_{K_{x,F}}$ and  $(\;,\; ;\;]_{K_{x,F}}$ :
\begin{equation} \label{diez}
\prod_{x \in F}(\phi_1,\phi_2,\phi_3)_{K_{x,F}}=1
\qquad \mbox{and} \qquad
\sum_{x \in F}(\phi_1,\phi_2;\zeta]_{K_{x,F}}=0
\end{equation}
for any $\phi_1,\phi_2,\phi_3 \in K_{x,F}^*$,
$\zeta \in K_{x,F}$.
Besides, in this product (this sum)
the number of terms $\ne 1$ ($\ne 0$) is finite.

Hence by means of expressions~(\ref{d}) and (\ref{dd})
we obtain that for fixed   $\xi \in K_s^*$
there exists only a finite number of points $x \in F$
such that  $(\fcf{\ph}{\ps}, \xi)_{K_s} \ne 1$
and      $(\fcf{\ph}{\ps}, \xi]_{K_s} \ne 0   \mbox{.}$

Now as in the proof of corollary~\ref{c1}
if we consider   $\xi = a \in k^*$, $a \ne 1$,
then we receive,
that there exists only a finite number of points  $x \in F$
such that $\nu_{K_s}(\fcf{\ph}{\ps}) \ne 0$.
If we consider  $\xi = \tau$,
then we receive,
that there exists a finite number of  points $x \in F$
such that $\nu_{K_s}(\fcf{\ph}{\ps}) = 0$,
but   $\fcf{\ph}{\ps} \notin {\cal U}^1_{K_s}$.
If we consider $\xi = \tau^{-n}$, $n \ge 1$,
then we obtain,
that for every $n \ge 1$
there exists only a finite number of points $x \in F$
such that
$\fcf{\ph}{\ps} \in {\cal U}^1_{K_s}$ and
$\nu_{K_s}(\fcf{\ph}{\ps}-1) = n$. Hence it follows, that
the product from corollary~\ref{c3}  converges
in $k((\tau))^*$.

Now the equality $ \prod\limits_{x \in F} \fcf{\ph}{\ps}=1  $
follows from~(\ref{diez}) and corollary~\ref{c1}.

\begin{cons}[The reciprocity law around a point.]   \label{c4}
Fix a point  $x \in X$,
 $\ph, \ps \in K_x^*$. Then
in the following product
the number of terms $\ne 1$ is finite and
\begin{equation}  \label{ddiez}
 \prod_{C \ni x} \fcc{\ph}{\ps}=1 \mbox{,}
\end{equation}
where $C$ runs all irredicuble curves   $C \subset X$
such that $x \in C$.
\end{cons}
{\bf Proof} (of corollary~\ref{c4}).

The equality~(\ref{ddiez}) is received by  the same methods
as analogous equality of corollary~\ref{c3}.
The finitness of terms of equality~(\ref{ddiez})
will follow from an explicit construction of $\fcc{\ph}{\ps}$,
when $C \ne F$. (See theorem~\ref{t2}).

Now we shall formulate the following theorem.

\begin{th}[ Explicit formulas.]  \label{t2}
The map $\fccn$, which is uniquely determined by corollary~\ref{c1}
of theorem~\ref{t1},
is given by the following explicit formulas.
\begin{enumerate}
\item  If $C \ne F$, then for any $\ph$ and $\ps $ from $K_{x,C}^*$
$$ \fcc{\ph}{\ps} \eqdef \nm\nolimits_{k(C)_x/K_s} \:
(\ph,\ps)_{K_{x,C}}  \mbox{,}$$
where $(\;,\;)_{K_{x,C}}$ is the tame symbol, which is
determined by the discrete valuation corresponding
to the local parameter $t_C$ of the field $K_{x,C}$.
\item If $C=F$, then for any $\ph$  and $\ps$ from $K_{x,F}^*$
$$ \fcf{\ph}{\ps} \eqdef \nm\nolimits_{k'((t))/k((t))} \:
((\ph,\ps)_{f,F})   \mbox{,}  $$
where a bimultiplicative map
 $(\ph,\ps)_{f,F} \; : \; K^*_{x,F} \times K^*_{x,F}
\to k'((t))^*$ is given by the following table ($K_{x,F} \simeq k'((u))((t))$).
\renewcommand{\arraystretch}{0}
\begin{equation} \label{tab1}
\begin{array}{|c||c|c|c|c|}
\hline
&\rule[-8pt]{0pt}{25pt}   u&  t&   b \in k'^* &
  \veps_2 \in {\cal E}^1_{K_{x,F}} \\
\hline \rule{0pt}{2pt} &  &&&  \\ \hline
\rule[-8pt]{0pt}{25pt}   u & -1 &t &b&
\exp  \mathop{\rm res}_u (\ln \veps_2 \frac{du}{u})  \\
\hline
\rule[-8pt]{0pt}{25pt}    t & t^{-1} & 1 &1& 1\\
\hline
\rule[-8pt]{0pt}{25pt}     a \in k'^* & a^{-1} & 1& 1 &1\\  \hline
\rule[-8pt]{0pt}{25pt}    \veps_1 \in {\cal E }^1_{K_{x,F}}&
\left(\exp  \mathop{\rm res}_u (\ln \veps_1 \frac{du}{u})\right)^{-1} &1 &1&
\left(\exp  \mathop{\rm res}_u (\ln \veps_1 \frac{d_u \veps_2}{\veps_2})\right)^{-1} \\
\hline
\end{array}
\end{equation}
\renewcommand{\arraystretch}{1}
Here the first column is the first argument of the map
$(\;,\;)_{f,F}$, the first row is the second argument of this map.
\end{enumerate}

By the multiplicative property and decomposition~(\ref{dec2}),
the map $(\;,\;)_{f,F}$
is extended to the whole group
$K^*_{x,F} \times K^*_{x,F}$.
\end{th}
Now if we verify conditions~1 and 2 of theorem~\ref{t1}
for the constructed map $\fccn$ from theorem~\ref{t2},
then we shall simultaneously prove theorems~\ref{t1} and \ref{t2}.
\begin{nt} {\em
Here $\mathop{\rm res}_u$ is a slightly modified map $\mathop{\rm res}_u$
from
section~\ref{ss1.1};
"new" $\mathop{\rm res}_u$  is a map from  $K_{x,F} \cdot d_u(K_{x,F})$
to $k'((t))$,
where for any $\vph$ from $K_{x,F}$
${d_u \vph= \frac{\partial \vph}{\partial u}\; du}$ :
 $$\mathop{\rm res}\nolimits_u (\sum_i u^i \cdot f_i(t) du) \eqdef f_{-1}(t)
 \in k'((t)) \quad \mbox{.}$$
(Compare with definition~\ref{ado2}).
All properties of the "old" map $\mathop{\rm res}_u$ (definition~\ref{ado2})
transfered without changing to the "new" map $\mathop{\rm res}_u$.
} \end{nt}
\begin{nt} {\em
In the definition of the map $f_*(\;,\;)_{x,F}$
we use  local parameters $u$ and $t$  in the table~\ref{tab1}.
It arises a question about the dependence
of this definition on the choice of the local parameters.
But $f_*(\;,\;)_{x,F}$   satisfy the conditions of theorem~\ref{t1}.
Therefore by corollary~\ref{c1} of theorem~\ref{t1}
this map does not depend on the choice of the local parameters
$u \in \hat{\cal O}_{x,X}$  and  $\tau \in K_s$.
(Remind that we fix the canonical embedding
  $f^* : K_s \hookrightarrow K_{x,F}$ and $t = F^*(\tau) $.)
}
\end{nt}
\begin{nt} {\em
From corollary~\ref{cor2} of theorem~\ref{t1},
corollary~\ref{c4} of theorem~\ref{t1}
and from table~(\ref{tab1})
it follows easily that the map $(\;,\;)_{f,F}$
from theorem~\ref{t2} satisfy the conditions of
lemma~7 from~\cite{K}.
Therefore the map $(\;,\;)_{f,F}$  is the back map of the
Kato's "residue homomorphism" introduced by him in~\cite{K}.
}
\end{nt}

\begin{nt} {\em
If we change in the table~(\ref{tab1}) all "letters" $u$
by "letters" $t$  and  all "letters" $t$
by "letters" $u$  we obtain the back map to the usualy
tame symbol.
}
\end{nt}

\begin{nt} {\em
Here is some properties of the table~(\ref{tab1}),
which will be useful later:
\begin{itemize}
\item This table is the skew-symmetric table
with respect to the diagonal.
\item The diagonal elements $(u,u)_{f,F}$ and $(t,t)_{f,F}$
is regenerated from the nondiagonal elements by means of
the equality $(\ph,\ph)_{f,F}=(-1,\ph)_{f,F}$.
\item The last row can be expressed by means of uniform
formula: for any $\ph \in {\cal E}^1_{K_{x,F}}$ and
$\ps \in K^*_{x,F}$
$$
(\ph,\ps)_{f,F}=\left(\exp  \mathop{\rm res}\nolimits_u
(\ln \ph \frac{d_u\ps}{\ps})\right)^{-1}  \mbox{.}
$$
\end{itemize}
} \end{nt}

In table~(\ref{tab1}) we use maps $\exp$ and $\ln$ in 2-dimensional
local fields. These maps are well defined by the following two
lemms.
\begin{lemma}
Let $K$ be a 2-dimensional local field and $m = \ee^1_K -1$. Then
the maps
$\ln : \ee^1_K \to m$ and $\exp : m \to \ee^1_K$
defined by the corresponding series converges in the topology of
2-dimensional local field. (See definition of this topology
in~\cite{PF} or \cite{P3}).
Besides, these maps are mutual isomorphismes  between  $\ee^1_K$
and $m$.
\end{lemma}
\proof
From~\cite{PF} (or \cite{P3}) for the convergence
of the series   $\ln$ and $\exp$ it is enough to prove that if
$x \in m$, then
${\vphantom{(}  x^n \stackrel{n \to \infty}{\longrightarrow} 0}$.
But the last is trivial. (See~\cite{PF} or \cite{P3}).

Now by the same method as in~\cite{BSh} we obtain
that the maps     $\ln$ and $\exp$ are mutual isomorphisms.

From the following lemma we obtain that the exponential map
in 1-dimensional local field from the table~(\ref{tab1})
is well defined.
\begin{lemma}
Let $K=k((u))((t))$ be a 2-dimensional local field,
$n$ is the maximal ideal of local ring $k[[u]]$.
Then for any $ \ph \in \ee^1_K$ and $\ps \in K$ the values of the expression
$\mathop{\rm res}_u(\ln \ph \frac{d_u \ps}{\ps})$ are in $n$.
\end{lemma}
\proof
The map $\mathop{\rm res}_u(\ln \ph \frac{d_u \ps}{\ps})$ is the
multiplicative map with respect to the second argument.
Therefore the proof will follow from the checking
on the multiplicative generators of the field $K$.
(See decomposition~(\ref{dec2}).)
This checking is by direct calculations and we omit it here.

Now we start to prove properties~(\ref{d}) and (\ref{dd}) for the
maps from theorem~\ref{t2}. For this goal
we shall need in the following well-known lemma.

\begin{lemma}   \label{l1}
Let $k((\tau)) \hookrightarrow k'((\tau'))$
be an embedding of fields.
Then for any $\zt$ from $k'((\tau'))^*$ and $\xi$ from $k((\tau))^*$
we have:
$$ \nm\nolimits_{k'/k}(\zt,\xi)_{k'((\tau'))} =
(\nm\nolimits_{k'((\tau'))/k((\tau))} \zt,
\xi)_{k((\tau))}  \quad  \mbox{,}
$$
where $(\;,\;)_{k'((\tau'))}$ and $(\;,\;)_{k((\tau))}$ are
the tame symbols in the 1-dimensional local
fields  $k'((\tau'))$ ¨ $k((\tau))$.
\end{lemma}
\proof
See~\cite[ ch.~3, \S4, lemma~3]{S}.\\[5pt]
{\bf Proof} of the property~(\ref{d}) from theorem~\ref{t1}
for the map $\fccn$  from theorem~\ref{t2}:
\begin{enumerate}
\item  If $C \ne F$,
then from lemma~\ref{l1} and definition of $\fccn$
we obtain that it will be enough to prove:
	$$ (\ph,\ps,f^*(\xi))_{K_{x,C}}=((\ph,\ps)_C,\xi)_{k(C)_x}$$
But the last expression is an obvious corollary of expression~(\ref{tc2})
and of the following expression:
$$ (\alpha,\beta,\gamma)_{K_{x,C}}=
(\partial_2(\alpha,\beta),\gamma)_{k(C)_x} \quad \mbox{,}
$$
where $\alpha$ and $\beta$ are from $K_{x,C}^*$, $\gamma$
is from $k(C)_x^*$.
\item If $C=F$,
then due to lemma~\ref{l1} it will be enough to verify
the following expression:
\begin{equation} \label{lf}
(\ph,\ps,\xi)_{K_{x,F}} = ((\ph,\ps)_{f,F},\xi))_{k'((t))}
\end{equation}
for any $\ph$, $\ps$ from $K_{x,F}^*$ and $\xi$ from $k'((t))^*$.
But from the multiplicative property of expression~(\ref{lf})
this expression will follow after its checking on
the multiplicative generators of the fields $K_{x,F}$
and $k'((t))$.
(See decompositions~(\ref{dec1}) and (\ref{dec2}).)
We have some cases .
But all cases appeared here is proved by easy
direct calculations and we omit it here.
The proof of property~(\ref{d}) is finished.
\end{enumerate}
$\bf Proof$ of the propertry~(\ref{dd}) from theorem~\ref{t1}
for the map $\fccn$ from theorem~\ref{t2}.

At first, it is an easy remark,
that the property~(\ref{dd}) is equivalent to the commutativity
of the following diagram:
\begin{equation}   \label{ff}
\begin{array}{ccc}
K_{x,C}^* \times K_{x,C}^* &
\begin{picture}(60,8)
 \put(0,2){\vector(1,0){60}}
\put(0,8){$\scriptstyle
(\ph,\ps) \mapsto \frac{d\ph}{\ph} \wedge \frac{d\ps}{\ps} $}
 \end{picture}
& \om^2_{K_{x,C}/k'}\\
 \begin{picture}(0,36)
 \put(0,33){\vector(0,-1){33}}
 \put(0,15){$\; \fccn $}
 \end{picture}
 &&
 \begin{picture}(0,36)
 \put(0,32){\vector(0,-1){32}}
 \put(0,14){$\; f_*^{x,C} $}
 \end{picture}\\
  K_s^* &
\begin{picture}(60,8)
 \put(0,2){\vector(1,0){60}}
\put(18,8){$
\scriptstyle
h \mapsto \frac{d h}{h}  $}
 \end{picture}
 &
 \om^1_{K_s/k}
 \end{array}
\end{equation}

At second, we shall need in the following well knowm lemma:
\begin{lemma}      \label{l2}
Let $k((\tau)) \hookrightarrow k'((\tau'))$
be an embedding of fields,
then the following diagram is commutative:
$$
\begin{array}{ccc}
k'((\tau'))^* & \stackrel {f \mapsto \frac{df}{f} }{\longrightarrow}
& \om^1_{k'((\tau'))/k'}\\
 \begin{picture}(0,36)
 \put(0,33){\vector(0,-1){33}}
 \put(0,15){$\; \nm $}
 \end{picture}
 &&
 \begin{picture}(0,36)
 \put(0,32){\vector(0,-1){32}}
 \put(0,14){$\; \tr $}
 \end{picture}\\
 k((\tau))^*  & \stackrel {g \mapsto \frac{dg}{g} }{\longrightarrow}     &
 \om^1_{k((\tau))/k}
\end{array}
$$
\end{lemma}

Now from lemma~\ref{l2},
diagram~(\ref{ff}) and remark~\ref{n5}
we obtain easily that property~(\ref{dd}) follows from the
following two lemms:
\begin{enumerate}
\item $C \ne F$
\begin{lemma} \label{l4}
Let $K=k((t))((u))$ be a 2-dimensional local field.
Then for any $\xi$ from $\bar K=k((t))$, $f$ and $g$ from $K^*$
we have
\begin{equation} \label{*}
\mathop{\rm res}\nolimits_{\bar K/k}\left(\xi \cdot \frac{d(f,g)_K}{(f,g)_K}\right)=
\mathop{\rm res}\nolimits_{K/k}\left(\xi \cdot \frac{df}{f} \wedge \frac{dg}{g}\right)
\quad \mbox{.}
\end{equation}
\end{lemma}
\item $C = F$
\begin{lemma}   \label{ll}
Let $K=k((u))((t))$ be a 2-dimensional local field.
Then for any  $\ph$ and $\ps$ from $K^*$, $\xi$ from $k((t))$
we have
\begin{equation} \label{f*}
\mathop{\rm res}\nolimits_{K/k}\left(\xi \frac{d \ph}{\ph} \wedge \frac{d \ps}{\ps}
\right)= \mathop{\rm res}\nolimits_{k((t))/k}
\left( \xi \frac{d (\ph,\ps)_{f,F}}{(\ph,\ps)_{f,F}}    \right)\quad \mbox{.}
\end{equation}
\end{lemma}
\end{enumerate}
{\bf Proof}  (of lemma~\ref{l4}).

Let us remark,
that the left and right hand  sides of formula~(\ref{*})
are additive expressions with respect to
 $\xi$
and multiplicative expressions with respect to $f$ and $g$.
Besides, left and right sides of formula~(\ref{*})
are symbol expressions,
i.~e. if $f=\ph$, $g=1-\ph$,
then these expressions are equal to $0$.
Therefore we shall consider cases when $\nu_t(\xi) \ge 0$
or $\xi = t^{-l}$ for all $l \ge 1$.

If $\nu_t(\xi) \ge 0$, then formula~(\ref{*}) is
expression~(6) from~\cite[ ch. 2.5.3]{PF} (or see~\cite[ prop. 2(5)]{P2},
\cite{L}).

Now let $\xi = t^{-l} $, $ l \ge 1 $. Then by the multiplicative
property it is enough to consider cases
when $f \in k^*$,
 $f \in \ee^1_K$, $f=t$, $f=u$ and
$g \in k^*$,
 $g \in \ee^1_K$, $g=t$, $g=u$.
By the skew-symmetric property with respect
to $f$ and $g$ some cases can be omitted.
All cases appeared here are proved by easy direct calculations.
Therefore we omit their proofs. As an example,
we consider only one case: let $f \in \ee^1_K$, $g=u$.
We have  $\ee^1_K = (1+ t \cdot k[[t]] )    (1+u \cdot k((t))[[u]])$,
therefore there exists a decomposition $f = e_1 \cdot e_2$,
where $e_1 \in 1+t \cdot k[[t]]$, $e_2 \in 1+u \cdot k((t))[[u]] $.
Then on the left hand side of~(\ref{*}) we have:
\begin{equation}  \label{cc}
\mathop{\rm res}\nolimits_{\bar K/k}   \left(t^{-l} \cdot
  \frac{\frac{\partial }{\partial t} (e_1 \cdot e_2 , u)_K  }
  {(e_1 \cdot e_2 ,u)_K } \: dt \right)=
  \mathop{\rm res}\nolimits_{\bar K/k}  \left(t^{-l}
  \cdot \frac{\frac{\partial}{\partial t} e_1}{e_1}    \:
   dt \right)  \quad \mbox{.}
\end{equation}
On the right hand side of~(\ref{*}) we have:
$$
\begin{array}{rl}
&
\mathop{\rm res}\nolimits_{K/k}\left(t^{-l} \cdot \frac{d(e_1 \cdot e_2)}{e_1 \cdot e_2}
\wedge  \frac{du}{u}\right)   =\\   =&
\mathop{\rm res}\nolimits_{K/k}\left( t^{-l} \cdot \frac{de_1}{e_1} \wedge
\frac{du}{u} \right)
+ \mathop{\rm res}\nolimits_{K/k}\left( t^{-l} \cdot \frac{de_2}{e_2}
\wedge \frac{du}{u}) \right)
=\\
=& \mathop{\rm res}_{\bar K/k} (t^{-l} \frac{\pat e_1}{e_1} dt )+
\mathop{\rm res}\nolimits_{K/k}( t^{-l} \cdot \frac{de_2}{e_2} \wedge \frac{du}{u})
\quad \mbox{.}\\
\end{array}
$$
Comparing the last expression with~(\ref{cc})
we obtain,
that it is enough to prove
$\mathop{\rm res}\nolimits_{K/k}( t^{-l}
\cdot \frac{de_2}{e_2} \wedge \frac{du}{u})=0$.
But
$t^{-l} \cdot \frac{de_2}{e_2} \wedge \frac{du}{u} =
t^{-l} \cdot \frac{\pat e_2}{e_2 \cdot u}  \wedge du$,
and $\frac{\pat e_2}{e_2} \in u \cdot k((t))[[u]]$,
therefore $\frac{\pat e_2}{e_2 \cdot u} \in  k((t))[[u]]$,
and $\mathop{\rm res}_{K/\bar K}(t^{-l} \cdot \frac{de_2}{e_2} \wedge
\frac{du}{u})=0$. The proof of lemma~\ref{l4} is finished.\\[12pt]
{\bf Proof} (of lemma~\ref{ll}).

As above,
from the additive property with respect to $\xi$,
the multiplicative property with respect to $\ph$ and $\ps$,
and the skew-symmetric property
with respect to $\ph$ and $\ps$ we have some cases:
\begin{enumerate}
\item $\nu_{k((t))} (\xi) \ge 0$, then the left part of~(\ref{f*})
is equal to
$ \quad \bar{\xi} \cdot {\rm det}
$$
\left (
\begin{array}{cc}
\nu_{\bar K}(\pi(\ph)) &   \nu_{ K}(\ph) \\
\nu_{\bar K}(\pi(\ps)) &    \nu_{ K}(\ps) \\
\end{array}
\right )
$$   $,
 where $\bar \xi \in k$, $\bar \xi = \xi \bmod (t \cdot k[[t]])$.
(See~\cite[ ch. 2.5.3, prop. 6]{PF}, or~\cite[ prop. 2(5)]{P2} )
We have several subcases:
\begin{enumerate}
\item Let $\nu_{k((t))}(\xi) > 0$, then $\bar \xi =0$ and
the left part of~(\ref{f*}) is equal to $0$.
By checking on multiplicative generators of $K$
we conclude,
that the right part of~(\ref{f*})  is also equal to $0$.
\item Let $\nu_{k((t))} (\xi) =0$. Then we can consider
 $\xi \in {\cal U}^1_{k((t))}$. Therefore we have to prove\\
$$ {\rm det}
\left (
\begin{array}{cc}
\nu_{\bar K}(\pi(\ph)) &   \nu_{ K}(\ph) \\
\nu_{\bar K}(\pi(\ps)) &    \nu_{ K}(\ps) \\
\end{array}
\right )
  =\;\mathop{\rm res}\nolimits_{k((t))/k}
  \left (\xi \frac{d(\ph,\ps)_{f,F}}{(\ph,\ps)_{f,F}} \right )
\quad \mbox{,}
$$
where $\pi : {\cal O}_K \to \bar{K}$ is the canonical map.

The last formula is proved by direct checking
on the multiplicative generators of $K$.
\end{enumerate}
\item  $\xi = t^{-k} \quad , \quad k \ge 1$.
We have again some uncomplicated cases on values of  $\ph$ and $\ps$.
We omit almost all cases here and consider only one.
Let  $\ph \in \ee^1_K$,
$\ps \in K^*$. Then expression~(\ref{f*}) is equivalent to the
following expression:
\begin{equation}     \label{**}
\mathop{\rm res}\nolimits_{u} \frac{d \ph}{\ph} \wedge \frac{d \ps}{\ps}     \quad
\stackrel{?}{=}        \quad
\frac{d_t(\exp \mathop{\rm res}_u (- \ln \ph \cdot \frac{d_u \ps}{\ps}))}
{\exp \mathop{\rm res}_u (- \ln \ph \cdot \frac{d_u \ps}{\ps})}    \quad \mbox{.}
\end{equation}
Transform the right part of~(\ref{**}):
$$
\begin{array}{ccccc}
&
\frac{d_t(\exp \mathop{\rm res}_u (- \ln \ph \cdot \frac{d_u \ps}{\ps}))}
{\exp \mathop{\rm res}_u (- \ln \ph \cdot \frac{d_u \ps}{\ps})} &
= &
d_t(\ln (\exp  \mathop{\rm res}_u (- \ln \ph \cdot \frac{d_u \ph}{\ph})))&
=  \\
= &
d_t ( \mathop{\rm res}_u (- \ln \ph \cdot \frac{d_u \ph}{\ph}))  &
=&
\mathop{\rm res}_u (\pat (-\ln \ph \cdot \frac{\partial \ps}{\partial u} \cdot
\frac{1}{\ps}))  \; du \wedge dt & \mbox{.}    \\
\end{array}
$$
Now~(\ref{**}) is equivalent to the following expression:
$$
      \mathop{\rm res}\nolimits_u \frac{d \ph}{\ph} \wedge \frac{d \ps}{\ps}
\quad      \stackrel{?}{=}     \quad
 \mathop{\rm res}\nolimits_u (\pat (-\ln \ph \cdot \frac{\partial \ps}{\partial u} \cdot
      \frac{1}{\ps})) \; du \wedge dt \quad \mbox{.}
$$
We shall compare the differential forms from left and right
hand  sides of the last expression:
 $\frac{d \ph}{\ph} \wedge \frac{d \ps}{\ps}$
and    $  \pat (-\ln \ph \cdot \frac{\partial \ps}{\partial u} \cdot
      \frac{1}{\ps})  \; du \wedge dt \quad \mbox{.}$ \\
The first differential form is equal to
\begin{eqnarray}
  \frac{d \ph}{\ph} \wedge \frac{d \ps}{\ps} &= &
\frac{1}{\ph \ps} (\frac{\partial \ph}{\partial u}du +
\frac{\partial \ph}{\partial t} dt) \wedge
(\frac{\partial \ps}{\partial u}du + \frac{\partial \ps}{\partial t} dt)=
\nonumber \\
& = &
\frac{1}{\ph  \ps}(\frac{\partial \ph}{\partial u} \cdot
\frac{\partial \ps}{\partial t} - \frac{\partial \ph}{\partial t}
\cdot \frac{\partial \ps}{\partial u}) \; du \wedge dt \; \mbox{.}  \label{b}
\end{eqnarray}
The second differential form is equal to
\begin{eqnarray}
&&\pat (-\ln \ph \cdot \frac{\partial \ps}{\partial u} \cdot
      \frac{1}{\ps}) \; du \wedge dt = \nonumber \\&=&
  (-\frac{\partial \ln \ph}{\partial t} \cdot
  \frac{\partial \ph}{\partial u} \cdot \frac{1}{\ps}) \; du \wedge dt
  - \ln \ph \cdot \pat(\frac{\partial \ps}{\partial u} \cdot \frac{1}{\ps})
\; du \wedge dt =   \nonumber \\
& = & -\frac{1}{\ph \ps} \cdot \frac{\partial \ph}{\partial t}
\cdot \frac{\partial \ps}{\partial u } \; du \wedge dt
-\ln \ph \cdot \pat(\frac{\partial \ps}{\partial u} \cdot \frac{1}{\ps})\;
du \wedge dt  \quad \mbox{.} \label{bb}
\end{eqnarray}
Substracting the expression~(\ref{bb}) from~(\ref{b}) we obtain:
\begin{eqnarray*}
& & \frac{d \ph}{\ph} \wedge \frac{d \ps}{\ps}   -
  \pat (-\ln \ph \cdot \frac{\partial \ps}{\partial u} \cdot
      \frac{1}{\ps})  \; du \wedge dt =\\  &=&
\frac{1}{\ph  \ps} \cdot \frac{\partial \ph}{\partial u} \cdot
\frac{\partial \ps}{\partial t}  \; du \wedge dt +
\ln \ph \cdot \pat\;(\frac{\partial \ps}{\partial u} \cdot \frac{1}{\ps})\;
du \wedge dt =\\ &=&
\frac{\partial (\ln \ph)}{\partial u} \cdot \frac {1}{\ps}   \cdot
\frac{\partial \ps}{\partial t}\; du \wedge dt +
\ln \ph \cdot \pat \; (\frac{\partial \ps}{\partial u} \cdot \frac{1}{\ps})
\; du \wedge dt= \\ &=&
\frac{\partial (\ln \ph)}{\partial u}
 \cdot \frac {1}{\ps} \cdot
\frac{\partial \ps}{\partial t}\; du \wedge dt +
\ln \ph \cdot (\pat \; \frac {\partial }{\partial u} \ps) \cdot \frac{1}{\ps}
\; du \wedge dt +\\
&+&\ln \ph \cdot
\frac{\partial \ps}{\partial u} \cdot \pat (\frac{1}{\ps})\; du \wedge dt=\\
&=&
\frac{\partial (\ln \ph)}{\partial u}
 \cdot \frac {1}{\ps}  \cdot
\frac{\partial \ps}{\partial t}\; du \wedge dt +
\ln \ph  \cdot (\frac{\partial}{\partial u} \;
\pat \ps ) \cdot \frac{1}{\ps}
\; du \wedge dt +\\
&+&\ln \ph \cdot \frac{\partial \ps}{\partial u}
\cdot (- \frac{1}{\ps^2} \cdot \frac{\partial \ps}{\partial t})\; du \wedge
dt=\\
& = &
\frac{\partial (\ln \ph)}{\partial u}
 \cdot \frac {1}{\ps}  \cdot
\frac{\partial \ps}{\partial t}\; du \wedge dt +
\ln \ph  \cdot (\frac{\partial}{\partial u} \;
\pat \ps ) \cdot \frac{1}{\ps}
\; du \wedge dt + \\
&+&\ln \ph \cdot \frac{\partial \ps}{\partial t} \cdot
\frac{\partial}{\partial u} \;(\frac{1}{\ps})\; du \wedge dt =\\
&=&
\frac{\partial (\ln \ph)}{\partial u}
 \cdot \frac {1}{\ps}  \cdot
\frac{\partial \ps}{\partial t}\; du \wedge dt +
\ln \ph \cdot \frac{\partial}{\partial u}\;(\frac{1}{\ps} \cdot
\frac{\partial \ps}{\partial t}) \; du \wedge dt=\\
&=&
\frac{\partial}{\partial u}\;(\ln \ph \cdot \frac{1}{\ps} \cdot
\frac{\partial \ps}{\partial t})\;du \wedge dt \quad \mbox{.}
\end{eqnarray*}
But for any $\zt$ from $K$
$$\mathop{\rm res}\nolimits_u\left(\frac{\partial }{\partial u}\;
(\zt) \; du \wedge dt\right) = 0 \quad \mbox{.} $$
Therefore expressions~(\ref{b}) and~(\ref{bb}) gives us the same
residues with respect to the parameter $u$.
The proof of lemma~\ref{ll} is finished.
\end{enumerate}

Now the proof of  property~(\ref{dd}) from theorem~\ref{t1} is finished.
Thus we proved theorem~\ref{t1} and theorem~\ref{t2}.

\section{$K_2$-- adeles, cohomology of $K_2$-- functors and $K_2$--
Gysin morphism}
Retain all notations from section~\ref{s2} and subsection~\ref{ss1.2}.
Note also that in this section we shall not assume $char k =0$
except for specially marked cases.

Besides, if  $K_{x,C} \simeq k((t_1))((t_2))$  is a 2-dimensional
local field, which associated to the pair  $x \in C$;
then by  $\oo_{x,C}$  denote the ring $\oo_{K_{x,C}}$,
i.~e.  $\oo_{x,C} \simeq k((t_1))[[t_2]]$.
Remind also, that by $\oo_C$ we have denoted  the ring $\oo_{K_C}$,
where $K_C \simeq k(C)((t_C))$ and $C \subset X$ is
an irredicuble curve.

 By $(\;,\;)_C$ denote the tame symbol associated with the irredicuble
curve $C$. (Such curve $C$ give the discrete valuation of the fields $k(X)$
and  $K_C$.)

 Let $\hat{K}_{x,X}$ be the fraction field of the ring $\oo_{x,X}$.

 Let $\hat \oo_{x,X}(\infty C) \eqdef ( \hat \oo_{x,X})_{(t_C)} $
be the localization of the ring  $\oo_{x,X}$
along the ideal $(t_C)$ of the curve $C$.

If $A$ is a ring, then let $K_2(A)$ be the Quillen $K$-functor.
If $\gamma$ is an ideal of the ring $A$,
then $K_2(A,\gamma) \eqdef {\rm Ker}\:(K_2(A)
\stackrel{p}{\to} K_2(A/ \gamma)) $.

By Stein result (see~\cite{M})  we have an explicit description
of functor $K_2$ for some rings:
\begin{enumerate}
\item
Let $A$ be a local ring,
then $K_2(A)$ is the group, generated by symbols $(u,v)$,
i.~e. $u$ and $v$  are from $A^*$ and satisfy
the usual symbol properties:
the bimultiplicative property and $(\ph,1-\ph)=1$
for $\ph$ and $1-\ph$ from $A^*$.
\item
Let $A$ be a local ring, $\gamma$ be an ideal in maximal ideal of $A$.
Then the natural map $p : K_2(A)  \to K_2(A/\gamma)$
is the quotient map of symbols,
and $K_2(A,\gamma)$ as subgroup of $K_2(A)$
is generated by the following symbols: $(1+i,u)$, where $i \in \gamma$,
$u \in A^*$.
\end{enumerate}

We have the following well known propositions:
\begin{prop}    \label{pk}
Let $K$ be the complete field of discrete valuation
with residue field $\bar K$, then
$$
 K_2({\cal O}_K) \cong {\rm Ker} \; (K_2(  K) \stackrel{(\;,\;)_{ K}}
 {\longrightarrow} {\bar K}^*  )   \quad \mbox{.}
$$
\end{prop}
\begin{prop}  \label{pp}
$$
K_2(\hat {\cal O}_{x,X} ) \simeq {{\rm Ker}\:} \left(K_2 ( \hat K_{x,X})
\stackrel{\bigoplus (\;,\;)_C}{\longrightarrow}
\bigoplus_C k(C)^* \right)  \quad \mbox{,}
$$
where this direct sum is over all prime ideals $C$ of height $1$ of the ring
$\hat {\cal O}_{x,X}$ ($\simeq$~irredicuble curves $C$
in ${ \rm Spec} \, \hat {\cal O}_{x,X}$), $k(C)$ is the fraction
field of quotient ring $\hat {\cal O}_{x,X}$ by $C$,
$(\;,\;)_C$  is the tame symbol associated with $C$.
\end{prop}
{\bf Remark} about the proof of these propositions.

From above explicit descriptions of $K_2$ in the symbol language
these propositions can be proved by means of noncomplete
direct calculations with symbols.
As an example,
we give later a such proof of proposition~\ref{pp}.

Note also that these both propositions follow from Gersten resolution.
For example, for the ring $\hat {\cal O}_{x,X}$
this resolution is proved in~\cite[ \S~7, theorems~5.6 and~5.13]{Q}
and its beginning is following:
$$
 0  \longrightarrow K_2(\hat K_{x,X})
 \stackrel{\bigoplus (\;,\;)_C}{\longrightarrow}
 \bigoplus_{C \subset X} K_1(k(C))\longrightarrow \ldots
$$
And the first cohomology group of this complex coinsides with
$$H^0 (\,{\rm Spec} \,\hat {\cal O}_{x,X}\,,\,
{\cal K}_2({\rm Spec} \, \hat {\cal O}_{x,X})\,) \simeq
K_2(\hat {\cal O}_{x,X}) \mbox{.}$$
(Here
    ${\cal K}_2({\rm Spec} \, \hat {\cal O}_{x,X}) $  is
the sheaf on  ${\rm Spec} \, \hat {\cal O}_{x,X}$
associated with the presheaf $\left \{ U \longmapsto K_2(U) \right.$;
$U$ is open set in $\left. {\rm Spec} \, \hat {\cal O}_{x,X} \right \}$.)\\[5pt]
{\bf Proof} of proposition~\ref{pp}.
From the symbol representation of $K_2(\hat{\oo}_{x,X})$
it is clear, that
$K_2(\hat {\cal O}_{x,X} ) \subset {{\rm Ker}\:} \left(K_2 ( \hat K_{x,X})
\longrightarrow
\bigoplus\limits_C k(C)^* \right) $.
Therefore we have to prove an back inclusion.
Let $a=\prod\limits_i (f_i,g_i) \in K_2(\hat K_{x,X})$
such that $a_C=\prod\limits_i (f_i,g_i)_C=1$
for any irredicuble curve $C$ in ${\rm Spec} \,\hat {\cal O}_{x,X}$.
We shall prove that  $a=\prod\limits_j (h_j,d_j)$,
 where $h_j$ and $d_j$ are from  $\hat {\cal O}_{x,X}^*$,
i.~e. $a \in K_2(\hat {\cal O}_{x,X})$.

 The ring $\hat {\cal O}_{x,X}$ is the ring
with a unique decomposition on prime factors.
Therefore from the multiplicative property of symbols
we have the following decomposition
 $a= \prod\limits_k a_k$, where $a_k = (p_k,q_k)$ and
$q_k$, $p_k$ are prime or invertible elements in $\hat {\cal O}_{x,X}$.

Fix some prime element $t$ from a collection  $\{p_k\}$ or $\{ q_k\}$.
Denote by $b$ the product those $a_k$,
for which either $p_k = t$ or $q_k = t$.

We shall prove, that $b \in K_2(\hat {\cal O}_{x,X})$.
From the multiplicative property,
the skew-symmetric property and the identity $(t,t)=(-1,t)$
one can suppose that
$b= (\gamma \cdot t_1^{m_1} \cdot \ldots \cdot t_l^{m_l},t) \cdot h$,
where $\gamma \in \hat {\cal O}_{x,X}^*$, $t_i$ are
prime elements from $\hat {\cal O}_{x,X}$ and $t_i \ne t$ for all $i$,
$h \in K_2(\hat {\cal O}_{x,X})$.
It is clear from decomposition of $b$,
that $(a \cdot b^{-1})_{(t)}=1$.
Besides, $(a)_{(t)}=1$.  Therefore
$(b)_{(t)}=1$.
Hence,  $\gamma \cdot t_1^{m_1} \cdot \ldots \cdot t_l^{m_l}=1+t \cdot d$
for some  $d \in \hat {\cal O}_{x,X}$.
But $1+t \cdot d \in \hat {\cal O}_{x,X}^* $,
hence $m_i=0$ for all $i$   from the uniqueness of decomposition
on prime factors.
Therefore $b=(1+t \cdot e,t) \cdot h$ for some
 $e \in \hat {\cal O}_{x,X}$.
Now we have two cases: 1) if $e \in \hat {\cal O}_{x,X}^*$,
then $(1+t \cdot e,t)=(1+t \cdot e,t)(1 + t \cdot e , -t \cdot e)^{-1}=
 (1+t \cdot e, -e^{-1}) \in K_2(\hat {\cal O}_{x,X})$;
2) if $e \notin \hat {\cal O}_{x,X}^*$,
then $(1+t \cdot e, t)= ((1-t)(1+t \cdot e),t)=(1+ t \cdot (-1+e-t \cdot e),t)$,
now $-1+e-t \cdot e \in \hat {\cal O}_{x,X}^*$,
therefore it is the previous case.

We proved  $b \in K_2(\hat {\cal O}_{x,X})$,
therefore $b_C=1$ for all $C$.
Hence $(a \cdot b^{-1})_C=1$ for all $C$,
and the symbol decomposition of $a \cdot b^{-1}$ does not contain
the prime element $t$.
Now the proof is by induction on prime elements
from collections  $\{ p_k\}$ and $\{ q_k \}$.
The proof is finished.

Remind that $K_x = k(X) \cdot \hat{\oo}_{x,X}$.
\begin{prop} \label{add1}
The following complex is an exact complex in the middle term:
\begin{equation} \label{kadd}
 K_2(K_x)
\stackrel{\varphi}{\longrightarrow}
K_2(\hat K_{x,X})
\stackrel{\bigoplus(\;,\;)_{\tilde C}}{\longrightarrow}
\bigoplus_{\tilde C} k(\tilde C)^*  \mbox{,}
\end{equation}
where the sum is over all prime ideals $\tilde C$ of height $1$ of
the ring $\hat{\oo}_{x,X}$
such that $\tilde C$  does not divide any prime ideal of
the ring $\oo_{x,X}$ under the natural embedding $\oo_{x,X}
\hookrightarrow \hat{\oo}_{x,X} $,
and $(\;,\;)_{\tilde C}$  is the tame symbol associated with $\tilde C$
\end{prop}
{\bf Proof.\ }

It is clear, that $K_x$  is 1-dimensional domain of principal
ideals. Therefore $K_x$   is the Dedekind ring,
whose maximal ideals are in one-to-one correspondence
to prime ideals $\tilde C$  of height~$1$
of the ring $\hat{\oo}_{x,X}$ such that $\tilde C$
does not divide any prime ideal of the ring $\oo_{x,X}$.
Besides,
 $\hat K_{x,X}$ is the fraction field of the ring $K_x$.
Now sequence~(\ref{kadd}) is a part
of the long $K$-sequence  for any
Dedekind ring, see~\cite{Q}.
Proposition~\ref{add1}  is proved.

Let us {\em denote} $\check K_2(K_x) \stackrel{\rm def}{=}
{{\rm Im}\:} (K_2(K_x) \stackrel{{\varphi}}
{\longrightarrow} K_2(\hat K_{x,X}))$.
By means of proposition~\ref{add1} it is not hardly to
describe the group $\check K_2(K_x)$ in the symbol language:
\begin{prop} \label{add2}
$\check K_2(K_x)$ as subgroup in $K_2(\hat K_{x,X})$
is generated by symbols $(u,v)$, where $u \in K_x^*$ and $v \in K_x^*$.
\end{prop}
{\bf Proof.\ }

From sequence~(\ref{kadd}) it follows that
the symbols
 $(u,v)$ $u \in K_x^*$,  $v \in K_x^*$
belong to $\check K_2(K_x)$.
The back inclusion is proved from sequence~(\ref{kadd})
by means of direct calculations with symbols.
That is done in a way completely analogous to that in the symbol proof
of proposition~\ref{pp}: by the same steps and
  by induction on prime elements of the ring
$\hat{\oo}_{x,X}$, which do not divide any prime elements of
the ring $\oo_{x,X}$. (See proof of proposition~\ref{pp}.)

\begin{prop}  \label{p9}
\begin{enumerate}
\item The following complex is exact in the middle term:
\begin{equation} \label{K*}
K_2(\hat K_{x,X}) \stackrel{\bigoplus (\;,\;)_C}{\longrightarrow}
\bigoplus_C k(C)^*
\stackrel
{{\textstyle \mathop{\bigoplus}\limits_C \nu_x}}{\longrightarrow}
{\mbox{\dbl Z}}
 \quad \mbox{,}
\end{equation}
where
this sum is over all prime ideals $C$ of height $1$
of the ring $\hat{\oo}_{x,X}$.
$(\;,\;)_{ C}$  is the tame symbol associated with $ C$,
and $\nu_x$   is the discrete valuation of function field
on the curve $C$ at the point $x$.
\item The following complex is exact in the middle term:
\begin{equation}    \label{K**}
\check K_2(K_x)
 \stackrel{\bigoplus (\;,\;)_{\breve C}}{\longrightarrow}
\bigoplus_{\breve C} k(\breve C)^*
\stackrel
{{\textstyle \mathop{\bigoplus}\limits_{\breve C} \nu_x}}{\longrightarrow}
{\mbox{\dbl Z}} \quad \mbox{,}
\end{equation}
where the direct sum is over all prime ideals  $\breve C$
of height~$1$ of the ring $\hat{\oo}_{x,X}$
such that $\breve C$  divides some prime ideal of the ring $\oo_{x,X}$
under the natural embedding $\oo_{x,X} \hookrightarrow \hat{\oo}_{x,X}$
($\simeq$~$\breve C$ are irredicuble curves in ${\rm Spec} \, \hat{\oo}_{x,X}$
such that these curves are in preimage of irredicuble
curves of  ${\rm Spec} \, {\oo}_{x,X}$
under the map ${\rm Spec} \, \hat{\oo}_{x,X} \to {\rm Spec} \, {\oo}_{x,X} $.)
\end{enumerate}
\end{prop}
{\bf The first proof.\ }

The sequence~(\ref{K*})
is Gersten resolution for  cohomologies of  sheaf
${\cal K}_2({\rm Spec} \, \hat {\cal O}_{x,X})$ on
${\rm Spec} \,\hat {\cal O}_{x,X}$. (See~\cite[ \S~7, theorem~5.13]{Q}.)
But we have $H^1 (\,{\rm Spec} \,\hat {\cal O}_{x,X}\,,\, {\cal F})=0$
for any sheaf $\cal F$ on  ${\rm Spec} \,\hat {\cal O}_{x,X}$.
Therefore the exactness in the middle term of sequence~(\ref{K*})
follows from
 $H^1 (\,{\rm Spec} \,\hat {\cal O}_{x,X}\,,\,
{\cal K}_2({\rm Spec} \, \hat {\cal O}_{x,X})\,) =
0$.

The exactness in the middle term of sequence~(\ref{K**}) follows
now from sequence~(\ref{K*}) and proposition~\ref{add1}.\\[10pt]
{\bf The second proof.\ }

That is an elementary proof in symbol language,
which don't use the methods of higher $K$-theory.

It is enough to prove~(\ref{K*}).
From~\cite[ ch. 6.1]{PF}
          we have  an elementary proof by means of blowing ups, that
\begin{equation} \label{c}
\left(\bigoplus\limits_C \nu_x \right) \circ
\biggl(\bigoplus (\;,\;)_C \biggr)=0
\end{equation}
Now  we shall prove the following fact:    for the collection
$\{f\} = \{f_C \::\: f_C \in k(C)^*\,;\, f_C=1,
\mbox{except a finite number irredicuble  curves $C$ in ${\rm Spec} \, \hat{\oo}_{x,X}$}\,;
\linebreak
\bigoplus\limits_C \nu_x (f_C)=0 \}$
there exist
 $g,h \in \hat K_{x,X}^*$ such that for any $ C $
$(g,h)_C=f_C$.
We shall prove it in two steps.\\[6pt]
\underline{Step~1.}

Suppose that all irredicuble curves $C$ in $\Soo$, such that $f_C \ne
1$ are regular curves at the point $x$.

Let $t_C$ be a local parameter at the point $x$ on curve
$ C $ in
${\rm Spec} \, \hat{\oo}_{x,X}$.
The proof will be by induction on a number of irredicuble
curves $ C $ in ${\rm Spec} \, \hat{\oo}_{x,X}$ such that   $f_C \ne 1$.
\begin{enumerate}
\item Let such $C$ be only one.
Choose an element $u \in \hat {\cal O}_{x,X}$
such that elements $t_C$ and $u$ are a pair of local parameters
at the point $x$, i.~e. $\hat {\cal O}_{x,X} \simeq k[[t_C,u]]$.
Let $k(C) \simeq k((t))$,
and let $i \::\: k(C) \hookrightarrow \hat K_{x,X}$
be an inclusion induced by the map
$t \stackrel{i}{\longmapsto} u$.
Then a symbol $(i(f_C),t_C)$ is a required symbol.
Indeed, we have
$$
\begin{array}{rcl}
(i(f_C),t_C)_C & = & f_C \\
(i(f_C),t_C)_{(u)}& =&1 \quad \mbox{.}
\end{array}
$$
We have $\nu_x(f_C)=0$, therefore
$i(f_C)$ is an invertible element in  $k[[u]]$. Hence   $i(f_C)$
is an invertible element in $\hat {\cal O}_{x,X}$.
Therefore we have
$$
\begin{array}{rcl}
(i(f_C),t_C)_E & =& 1 \quad \mbox{,}
\end{array}
$$
where $E$ is any other curve, $E \ne C$,
$E \ne (u)$.
\item
Suppose we have a collection from $n$ irredicubles curves
$C_i \,(i=1,\ldots, n)\; :\; f_C \ne 1 $.
Let $\Psi \in K_2(\hat K_{x,X})$ be
those unknown symbol,
which we have to construct. We shall calculate it.

Suppose that there exist pair $i \ne j$
such that the tangent lines of the curves $C_i$ and $C_j$
don't coinside,
i.~e. $C_i$ and $C_j$ are transversal curves in
${\rm Spec} \, \hat{\oo}_{x,X}$.
Without loss of generality  it can be assumed
that $i=1$, $j=2$.
Then $t_{C_1}$ and $t_{C_2}$
are the pair of local parameters
and
 $\hat {\cal O}_{x,X} \simeq k[[t_{C_1},t_{C_2}]]$.
Let $k(C_i) \simeq  k((t_i))$, and let
$i_{C_1}\: :\: k(C_1)\hookrightarrow \hat K_{x,X}$--
be an inclusion induced by the map
$t_1 \stackrel{i_{C_1}}{\longmapsto } t_{C_2}$.
Now consider a symbol $\Theta_1=
(i_{C_1}(f_{C_1}),t_{C_1}) \in K_2(\hat K_{x,X})$.
It is clear, that
$$
\begin{array}{rcl}
(i_{C_1}(f_{C_1}),t_{C_1})_{C_1}&=& f_C \\
(i_{C_1}(f_{C_1}),t_{C_1})_E&=& 1 \quad \mbox{,}
\end{array}
$$
for any curves $E \ne C_1$, $E \ne C_2$.
$\Theta_1 \in K_2(\hat K_{x,X})$,
therefore this symbol satisfy property~(\ref{c}).
$\Psi$ satisfy property~(\ref{c}) by condition.
Therefore  $\Theta_1^{-1} \Psi$
satisfy the property~(\ref{c}) and
$$
\begin{array}{rcl}
(\Theta_1^{-1} \Psi)_{C_1}&=&1\\
(\Theta_1^{-1} \Psi)_E & =& 1 \quad
\end{array}
$$
for any curves $E \ne C_2, \ldots , C_n$.
Now we apply the inductive hypothesis
to the $n-1$ curves $C_2, \ldots, C_n$
and collection
$$
f_{C_2}= (\Theta_1^{-1} \Psi)_{C_2}\; , \ldots , \;
f_{C_n}= (\Theta_1^{-1} \Psi)_{C_n}     \quad \mbox{.}
  $$
By inductive hypothesis
there exists $ \Theta_2  \in K_2(\hat{K}_{x,X})$
and we can write
$$
\begin{array}{rcl}
\Theta_1^{-1} \Psi & =& \Theta_2 \quad \mbox{, hence}\\
   \Psi & =&\Theta_1 \Theta_2
\end{array}
$$
And in this case~(\ref{K*}) is proved.

Now suppose that all curves $C_i \: (i=1, \ldots , n)$
have the same tangent line.
Then we consider the curve $G$ in ${\rm Spec} \, \hat{\oo}_{x,X}$,
transversally intersected to the all curves
$C_i\: (i=1, \ldots , n)$.
Now after  operations with the  pair
$C_1$ and $G$ such as above with  $C_1$ and $C_2$
we receive a case of the collection
$\{G,C_2, \ldots ,C_n \} $,
where the curve $G$ is transversal
to the all other curves.
But this case was yet considered.
Step~1 is proved.\\[6pt]
\end{enumerate}
\underline{Step~2.}

Now suppose that there exists irredicuble singular curve $C$
in $\Soo$, such that $f_C \ne 1$. We shall prove the following fact:
fix the smooth curve $C$, then there exist $f_1, f_2 \in  \hat K_{x,X}^*$,
such that  $(f_1, f_2)_C = f_C$.
And if $(f_1, f_2)_E \ne 1$ for some irredicuble curve $E$ in $\Soo$,
then $E$ is regular curve. Thus we reduce step~2  to step~1.

Fix some  regular local parameters $u$ and $t$,
then $\hat{\oo}_{x,X} \simeq k[[u,t]]$.
The ring  $k[[u,t]]$ is the ring with a unique decomposition on prime
factors. Also by Weierstrass preparation theorem (see~\cite[ ch.~VII, \S~3.8]{B})
we have, that every element $f \in k[[u,t]]$
has the following decomposition
$f = b \cdot u^{\nu_1(f)} \cdot g$,
where $g = t^{\nu_2(f)} + a_1 t^{\nu_2(f) -1} + \ldots + a_{\nu_2(f)}$,
$a_i  \in u \cdot k[[u]]$,
$b$ is an invertible element in $k[[u,t]]$,
$\nu_1(f) \ge 0$, $\nu_2(f) > 0$ (if $\nu_2(f) = 0$, then $g=1$).
Here $\nu_1(f)$ and $\nu_2(f)$  is uniqueli defined by $f$.
And $f \in k[[u,t]]$ is a prime element in this ring,
if and only if or $g=1, \nu_1(f)=1$,
or $\nu_1(f)=0$ and  $g$ is a prime polinom in the ring $k[[u]][t]$.
Also remark, that if $\nu_2(f) = 1$ and $\nu_1(f) = 0$,
then $f$ is a prime element in $k[[u,t]]$,
and $f$ is a regular curve in $\Soo$.

Let $t_C \in k[[u,t]]$ be a local parameter on curve $C$ in $\Soo$,
i.~e. $(t_C) = C$.
We shall do  induction on $\nu_2(t_C)$ (or $\nu_2(t_E)$).
If $\nu_2(t_C) =1$, then we have step~1 by the remark above.

Let   $\nu_2(t_C) > 1$.
Fix some elements $h_{1,C} , h_{2,C} \in k[[u,t]] $,
such that
$$
		 (h_{1,C} \, ({\rm mod} (t_C))) \cdot  (h_{2,C} \, ({\rm mod} (t_C)))^{-1} =f_C.
$$

If $\nu_2(h_{1,C}) < \nu_2(t_C)$ and
$\nu_2(h_{2,C}) < \nu_2(t_C)  $,
then in decomposition on prime elements we have
$$
h_{1,C} =\prod_i{p_i}
$$
$$
h_{2,C} =\prod_j{q_j}    \mbox{,}
$$
where for any $i,j$
$$
\nu_2(p_i) \le   \nu_2(h_{1,C})  <   \nu_2(t_C)
$$
$$
\nu_2(q_j) \le   \nu_2(h_{2,C})  <   \nu_2(t_C)  \mbox{.}
$$
Also for the symbol
$$
(h_{1,C} \cdot h_{2,C}^{-1} \:,\: t_C)
$$
we have, that
if  $E \ne (u)$, $E \ne C$,  $E \ne (p_i)$, $E \ne (q_j)$, then
$$
(h_{1,C} \cdot h_{2,C}^{-1} \:,\: t_C)_E =1   \mbox{.}
$$

Hence by inductive hypothesis  we conclude proof in this case.

 Otherwise, we can multiply $h_{i,C}$ by some degree of $u$,
and subtract from it an element $t_C$, multiplied by some invertible
element, by some degree of $u$, and by some degree of $t$.
And obtained elements $\tilde{h}_{i,C}$  have the following properties:
$$
		 (\tilde{h}_{1,C} \, ({\rm mod} (t_C))) \cdot  (\tilde{h}_{2,C} \, ({\rm mod} (t_C)))^{-1} =f_C.
$$
$$
 \nu_2(\tilde{h}_{1,C}) < \nu_2(t_C)
$$
$$
\nu_2(\tilde{h}_{2,C}) < \nu_2(t_C)
$$
But that is the case above.
We proved step~2 and proposition~\ref{p9}.\\[7pt]

Remind that if $x$ is a singular point of the irredicuble curve $C \subset
X$,
then
 $$K_{x,C} = \bigoplus_i K_{x,C_i}  $$
 $$\oo_{x,C}  = \bigoplus_i \oo_{x,C_i} \quad \mbox{,} $$
where $C_i$ are the "analitic" branches of the curve $C$ at the point $x$.
Therefore
$$
K_2(K_{x,C}) = \prod_i K_2(K_{x,C_i})
$$
$$
K_2(\oo_{x,C}) = \prod_i K_2(\oo_{x,C_i})
$$
$$
\nu_x (\;,\;)_C  \eqdef \sum_{C_i \ni x} \nu_x (\;,\;)_{C_i}  \quad \mbox{,}
$$
where $(\;,\;)_{C_i}$ is the usual tame symbol of the discrete valuation
field $K_{x,C_i}$, and $\nu_x$ is the discrete valuation on the branche $C_i$
at the point $x$.

\begin{defin}[$K_2$-adeles]   \label{kad}
Let $X$ be a surface and a pair $x \in C$ runs all irredicuble
curves $C \subset X$ and points $x \in C$, then
$$
\ka \eqdef \left \{ (\ldots,f_{x,C},\ldots) \in \prod_{x \in C}
K_2(K_{x,C})\; \mbox{ under the following conditions} \right\}
$$
\begin{enumerate}
\item
$f_{x,C} \in K_2(\oo_{x,C})$ for almost all curves $C \subset  X$.
\item
For all curves $C \subset X$, all integers $l \ge 1$ and all,
except a finite number (which depends on $l$) non singular points
$x \in  C$ we have
 $f_{x,C} \in K_2(\oo_{x,C}, m_C^l) \cdot
K_2(\hat \oo_{x,X}(\infty C))$,
where $m_C$ is the maximal ideal of the ring $\oo_{x,C}$.
\end{enumerate}
\end{defin}

\begin{nt} {\em
Definition 10 is similar to the definition from~\cite{DAN},
where it have been used for  global  class field theory on surface $X$.
See also~\cite{PF}} \end{nt}
\begin{prop}  \label{pp12}
Let $f=\{ f_{x,C}\}$ and $g=\{ g_{x,C}\}$ be from $\da^*_X$.
Then the element $(f,g)$
belongs to $\ka$, i.~e. there exists the natural computation
for
$K_2$-functor:
$$
\da_X^*
 \otimes  \da_X^* \longrightarrow \ka  \quad \mbox{.}
$$
\end{prop}
The proof of this proposition is enough easy and we omit it here.

\begin{prop}   \label{p11}
Let $(f,g) \in \ka$. Then under the Tate map:
$$
(f,g) \longmapsto \frac{df}{f} \wedge \frac{dg}{g}\in
\prod_{x \in C} \om^2_{K_{x,C}/k}
$$
the second adelic condition from definition~\ref{d8}
is true for the differential form
 $\frac{df}{f} \wedge \frac{dg}{g} $.
\end{prop}
\proof We have
$\frac{df}{f} \wedge \frac{dg}{g}=
\frac{1}{f\cdot g} (\frac{\partial f}{\partial u} \cdot
   \frac{\partial g}{\partial t} -
   \frac{\partial f}{\partial t} \cdot
   \frac{\partial g}{\partial u}
   )\, du \wedge dt   $,
 therefore
\begin{itemize}
\item if $(f,g) \in K_2(\hat \oo_{x,X}(\infty C))$,
then it is not hardly to prove that $\frac{df}{f} \wedge \frac{dg}{g}$
satisfies the needed condition.
\item if $(f,g) \in K_2(\oo_{x,C}, m_C^l)$ and
$\frac{df}{f} \wedge \frac{dg}{g}= \omega_{x,C}$,
then $\omega^{x,C}_{l-1}(u_{x,C})$ satisfies the needed condition.
(See definition~\ref{d8}).
\end{itemize}
The proposition~\ref{p11} is proved.
\begin{th}   \label{th3}
Using the diagonal map of $K_2(K_C)$ to $\prod\limits_{x \in C} K_2(K_{x,C})$
and of $K_2(K_x)$ to $\prod\limits_{C \ni x} K_2(K_{x,C})$ we put
$$
\begin{array}{rclrl}
K'_2(\da_1)&\eqdef& \prod\limits_{C \in X}K_2(K_C) \cap \ka &\subset & \ka\\
K'_2(\da_2) & \eqdef & \prod\limits_{x \in X}K_2(K_x)  \cap \ka &\subset & \ka\\
K'_2(\da_X(0)) & \eqdef & \prod\limits_{x \in C}K_2(\oo_{x,C}) \cap \ka &\subset & \ka\\
\kao & \eqdef & \prod\limits_{C \subset X}K_2(\oo_C) \;\cap \ka &\subset & \ka\\
\kaoo & \eqdef & \prod\limits_{x \in X}K_2(\hat \oo_{x,X} ) \;\cap \ka  &\subset & \ka
\quad \mbox{.}
\end{array}
$$
and consider the following complex $K_2({\cal A}_X)$:
{
\arraycolsep=2pt
$$
\begin{array}{@{}ccccc@{}}
K_2(k(X)) \times \kao \times
\kaoo & \longrightarrow  &
\kadd \times \kad \times \kado & \longrightarrow & \ka \\[7pt]
(f_0,f_1,f_2) & \mapsto  &
(f_2 \cdot f_0^{-1}, f_0 \cdot f_1, f_1^{-1}\cdot f_2^{-1}) \\
&& (g_1,g_2,g_3) & \mapsto & g_1 \cdot g_2 \cdot g_3
\end{array}
$$
}
Then $H^*(K_2({\cal A}_X)) \simeq  H^*(X,{\cal K}_2(X))$
and $H^2(K_2({\cal A}_X)) \simeq C\! H^2(X)$,
where   ${\cal K}_2(X) $  is  the sheaf on $X$
associated with the presheaf
$\left \{ U \longmapsto K_2(U) \right.$;
$U$ is open set in $\left. X \right \}$,
and $C\! H^2(X)$
is the  Chow group of algebraic cycles of degree $2$
with respect to the rational equivalence.
\end{th}
\begin{nt}    {\em
Now due theorem~\ref{th3} and proposition~\ref{pp12}
it can be possible to overwrite the adelic intersection
index of divisors on surface $X$ from~\cite{PF} (or~\cite{Lom} )
to the language of $K_2$-adeles.
}
\end{nt}
{\bf Proof} of theorem~\ref{th3}.

We have the following Gersten resolution on the surface $X$:
\begin{equation} \label{zv}
K_2(k(X)) \longrightarrow \bigoplus_{C \subset X}K_1(k(C)) \longrightarrow
\bigoplus_{x \in X}K_0(k(x)) \quad \mbox{.}
\end{equation}
And from~\cite[ \S~7, theorems~5.6 and~5.11]{Q} we have that $i$-th cohomology group
of this complex coinsides with
$H^i(X,{\cal K}_2(X))$.
We shall prove that this complex~(\ref{zv})
is quasi-isomorphic to complex  $K_2({\cal A}_X)$.
We consider the following diagram:
\begin{tabbing}
$K_2(k(X)) \times $\=$\kao \times$\=%
$\kaoo$\=$\longrightarrow $\=%
$\kadd \times $\=$\kad
\,\;\;\quad \quad\times$\=$\kado$\=$\longrightarrow $\=$\ka$\kill%
$K_2(k(X))$\>\>\>%
$\longrightarrow $\>\>$ \bigoplus\limits_{C \subset X}K_1(k(C))$%
\>\>$\longrightarrow$\>%
$\bigoplus\limits_{x \in X}K_0(k(x)) $\\%
\begin{picture}(16,36)
\put(16,1){\vector(0,1){39}}
\put(16,18){$\; {\ph_1 }$}
\end{picture}\>\>\>\>\>%
\begin{picture}(35,36)
\put(35,1){\vector(0,1){38}}
\put(35,18){$\; {\ph_2 }$}
\end{picture}\>\>\>%
\begin{picture}(34,36)
\put(34,1){\vector(0,1){38}}
\put(34,18){$\; {\ph_3 }$}
\end{picture}%
\\%
$K_2(k(X))\, \times $\>$ \kao \, \times $\>%
$\kaoo $\>$\longrightarrow$\>%
$\kadd \,\times $\>$\quad \;\,\kad\: \, \quad
\times $\>$\kado$\>$\longrightarrow $\>$\quad \: \ka $\\%
\>%
\begin{picture}(47,36)
\put(47,2){\vector(0,1){31}}
\put(47,16){$\; {\vph_1 }$}
\end{picture}\>\>\>\>%
\begin{picture}(35,36)
\put(35,2){\vector(0,1){31}}
\put(35,16){$\; {\vph_2 }$}
\end{picture}\>\>\>%
\begin{picture}(34,36)
\put(34,2){\vector(0,1){31}}
\put(34,16){$\; {\vph_3 }$}
\end{picture}%
\\%
\>$\kao \,\times$\>%
$\kaoo $\>$ \longrightarrow $\>%
$\kadd \,\times$\>$\quad\;\,\kao \:\, \quad
\times $\>$\kado$\>$\longrightarrow $\>$ \qquad \,\:\Phi $   ,\\%
\end{tabbing}
where $\Phi \eqdef \Ker \ph_3 \; \subset \ka$ and
maps $\ph_1 ,\ph_2 , \ph_3 , \vph_1 , \vph_2 , \vph_3$
are given by the following mode:
$$
\begin{array}{cccc}
\ph_1 \;:&
K_2(k(X))\, \times \, \kao \, \times \,
\kaoo & \longrightarrow & K_2(k(X))\\
&g_1\times g_2 \times g_3 & \mapsto & g_1 \\[20pt]
\ph_2\; :  &
\kadd \,\times \,\kad \,
\times \, \kado & \longrightarrow &
 \bigoplus\limits_{C \subset X}K_1(k(C))\\
& g_1 \times g_2 \times g_3 & \mapsto & \bigoplus\limits_{C \subset X}
(g_2)_{C} \\[20pt]
\ph_3 \;:& \ka & \longrightarrow& \bigoplus\limits_{x \in X}K_0(k(x))\\
& (\alpha_{x,C},\beta_{x,C}) &
   \mapsto &
\bigoplus\limits_{x \in X}(\sum\limits_{C \ni x}
\nu_x\:(\alpha_{x,C},\beta_{x,C})_C)\\
\end{array}
$$
$$
\begin{array}{cccc}
\vph_1 \; : &
 \kao \, \times \,
\kaoo & \longrightarrow &
K_2(k(X))\, \times \, \kao \, \times \,
\kaoo \\
& f_1 \times f_2 & \mapsto &
1 \times f_1 \times f_2 \\[20pt]
\vph_2 \; : &
\kadd \,\times \,\kao \,
\times \, \kado & \longrightarrow &
\kadd \,\times \,\kad \,
\times \, \kado \\
& g_1 \times g_2 \times g_3 &\mapsto& g_1 \times g_2 \times g_3 \\[20pt]
\vph_3 \; :
& \Phi & \longrightarrow & \ka\\
& \Ker \ph_3 &\hookrightarrow & \ka
\end{array}
$$
It is clear, that all three vertical complexes are exact sequences.
(The exactness of the middle complex follows from proposition~\ref{pk}).
Therefore the quasi-isomorphism
between the first horizontal complex and the second
horizontal complex will follow from the exactness
of the lower horizontal complex.
But the exactness of the lower horizontal
complex in the left term is obvious,
the exactness in the term $\Phi$
follows from sequence~(\ref{K**}) of proposition~\ref{p9},
the exactness in the middle term follows
from proposition~\ref{pp}.
Also  it is well known from the Bloch result,
that
$ H^2(X,{\cal K}_2(X)) \simeq C\! H^2(X)$.
This concludes the proof of theorem~\ref{th3}.\\[5pt]

In the sequel assume $char k = 0$.\nopagebreak

{\em Define} the map
 $f_* \eqdef
\mathop{ \prod\limits_{f(x)=s} }\limits_{C \ni x} \fccn$
from
 $\ka \subset \prod\limits_{x \in C}K_2(K_{x,C})  $
to the idele group $\da_S^* \subset \prod\limits_{s \in S} K_s^* $
of the curve $S$ for the projective morphism surface $X$ on the curve $S$.
(Compare with definition of the map $f_*$, see page~\pageref{adf}.)
\begin{prop}
The map
 $f_* \eqdef
\mathop{ \prod\limits_{f(x)=s} }\limits_{C \ni x} \fccn$
from
 $\ka   $
	to  $\da_S^*$
is well defined,
i.~e. this infinite product converges at every point $s \in S$.
\end{prop}
{\bf Proof.}

Over each point $s \in S$ we have
$$
 \mathop{\prod\limits_{f(x)=s}}\limits_{C \ni x}\fccn =
  \mathop{ \mathop{\prod\limits_{f(x)=s}}\limits_{ C \ni x}}\limits_{C \ne F}
  \fccn                       \times
    \mathop{\prod\limits_{f(x)=s}}\fcfn
  \quad \mbox{.}
$$

Now the first product from the right part of
the last expression contains only a finite
number terms $\ne 0$.
It follows from the first adelic condition of
definition~\ref{kad}.

The second product from the right part
of this expression converges.
That is proved by the same mode
as corollary~\ref{c3} from theorem~\ref{t1} by means
of proposition~\ref{p11} and
formulas~(\ref{d}) and (\ref{dd}) from theorem~\ref{t1}.
\begin{nt} {\em
The expression   $
\mathop{ \prod\limits_{f(x)=s} }\limits_{C \ni x} \fccn$
applied to the whole $\prod\limits_{x \in C} K_2(K_{x,C})$
make not sense,
since this series  will not converge.
}
\end{nt}
\begin{prop}  \label{qqu}
$f_*$ maps the complex $K_2({\cal A}_X)$ to the
complex
$$
1 \longrightarrow k(S)^* \times \prod_{s \in S } \hat{\oo}_s^*
\longrightarrow   \da_S^*
$$
and this map is a morphism of complexes.
\end{prop}
{\bf Corollary}
{\em
$f_*$ gives us the maps from
$H^1(X,{\cal K}_2(X))$
to $H^0(S,{\cal O}_S^*)$
and from             $ H^2(X,{\cal K}_2(X))$
to $H^1(S,{\cal O}_S^*)$.
}\\[5pt]
{\bf Proof} (of corollary)
It follows from proposition~\ref{qqu}  and theorem~\ref{th3}.\\[5pt]
{\bf Proof} (of proposition~\ref{qqu}).
It follows from corollary~\ref{c3}       and
corollary~\ref{c4}
of theorem~\ref{t1}.
Also, an explicit construction of map  $\fcfn$
(see table~(\ref{tab1}) from theorem~\ref{t2}) gives us that
$\fcfn \Bigl ( K_2(\oo_{x,F}) \Bigr )
\,\subset \prod\limits_{s \in S} \hat \oo_s^*$.
\begin{prop}  \label{ep}
The constructed map
 $f_* \: :\: H^2(X,{\cal K}_2(X)) \to H^1(S,{\cal K}_1(S))$
is the Gysin map from $C\!H^2(X)$ to $C\!H^1(S)$.
\end{prop}
\proof
At first, remind
that the Gysin map       $C\!H^2(X) \to C\!H^1(S)$
maps the point on the surface $X$
to the image point   on the curve $S$
with multiplicity equal to the degree of the extension
of the residue field of this point
over the residue field of the image point.

At second,
we have explicit isomorphismes
between the groups $H^2(X,{\cal K}_2(X))$
(in adelic variant)  and     $C\!H^2(X)$,
also between $H^1(S,{\cal K}_1(S))$ (in idelic variant)
and  $C\!H^{1}(S)$.
Moreover, we have the following diagram
$$
\begin{array}{@{}c@{}c@{}c@{}c@{}c@{}}
H^2(X,{\cal K}_2(X)) & \cong & \ka/\Bigl (\kadd \times \kad \times \kado \Bigr) &
\stackrel{\bigoplus\limits_{x \in X}
(\sum\limits_{C \ni x} \nu_x\, (\;,\;)_C)}{\longrightarrow}&
CH^2(X)\\ &&
 \begin{picture}(0,36)
 \put(0,33){\vector(0,-1){40}}
 \put(0,11){\llap {$\; {f_* }$}}
 \end{picture} &&
 \begin{picture}(0,36)
 \put(0,33){\vector(0,-1){40}}
 \put(0,18){\llap {$\; \mbox{Gysin's}$}}
 \put(0,7){\llap {$\; \mbox{map}$}}
 \end{picture}
 \\
H^1(S,{\cal K}_1(S)) & \cong &
\da^*_S /\Bigl (k(S)^* \times \prod\limits_{s \in S} \hat \oo^*_s \Bigr ) &
\stackrel{\bigoplus\limits_{s \in S} \nu_s }{\longrightarrow} &
C\!H^1(S)
\mbox{.}
\end{array}
$$
Now from the explicit construction
of the map $\fccn$ (see theorem~\ref{t2})
we conclude,
that this diagram is commutative diagram.
Therefore $f_*$ is the Gysin map.
The proposition is proved.


\nopagebreak
Moscow State University\\
{\em E-mail :} Dosipov@nw.math.msu.su

\begin{thebibliography}{99}
  \bibitem{BSh}
  Z.~I.~Borevich, I.~R.~Shapharevich  {\em Number theory.}
  Moscow, Nauka, 1985.
  \bibitem{B} N.~Bourbaki {\em Algebre Commutative}, Hermann, Paris, 1961-1965.
  \bibitem{K}
  K.~Kato {\em Residue homomorphism in Milnor K-theory}, Advanced studies
  in pure mathematics {\bf2} (1983) Galois groups and their representations,
  153-172.
    \bibitem{Lom}
        V.~G.~Lomadze {\em On the intersection index of divisors},
  Izv. Akad. Nauk SSSR Ser. Mat. 44 (1980),  1120-1130;
  English transl. in Math. USSR Izv. 17 (1981).
        \bibitem{L}
        V.~G.~Lomadze {\em On residues in algebraic geometry},
        Izv. Akad. Nauk SSSR Ser. Mat. 45(6) (1981),  1258-1287;
				English transl. in Math. USSR Izv. 19 (1982).
  \bibitem{MI}
  J.~Milnor {\em Algebraic K-theory and quadratic forms}, Invent. Math,
  9 (1969/70), 318-344.
  \bibitem{M}
  J.~Milnor {\em Introduction to algebraic K-theory}, Princeton
  Univ. Press, Princeton, N.~J., and Univ. of Tokyo Press, Tokyo, 1971.
\bibitem{UMN2} A.~N.~Parshin {\em Class fields and algebraic K-theory},
Uspekhi Matem. Nauk, v.30 (1975), 253-254.
\bibitem{P2}
        A.~N.~Parshin {\em On the arithmetic of two-dimensional
        schemes, I. Repartitions and residues}, Izv. Akad. Nauk SSSR Ser. Mat. 40(4) (1976), 736-773;
				English transl. in Math. USSR Izv. 10 (1976).
				\bibitem{DAN}
				A.~N.~Parshin {\em Abelian coverings of arithmetical schemes},
				DAN USSR, v.243 (1978), 855-858; English transl. in
				Soviet. Math. Doklady 19 (1978).
        \bibitem{P3}
        A.~N.~Parshin   {\em Local class field theory},
        Trud. Mat. Inst. AN SSSR, 165 (1984), 143-170;
				English transl. in Proceedings of the Steklov Institute
				of Mathematics 1985, Issue 3.
  \bibitem{PF}
  A. N.~Parshin, T. Fimmel {\em An introduction to the higher adelic
  theory},
  preprint.
  \bibitem{Q}
  D. Quillen, {\em Higher algebraic K-theory 1}, Lecture notes 341,
  Algebraic K-theory 1, 85-147.
  \bibitem{S}
  J.-P. Serre {\em Groupes algebriques et corps de classes},
  Hermann, 1959.
        \bibitem{Y}
        A. Yekutieli {\em An explicit construction of the Grothendieck
  residue complex}, Asterisque, 208 (1992).
\end{thebibliography}
\end{document}